\documentclass[final,3p,times]{elsarticle}

\usepackage{lineno,hyperref}
\modulolinenumbers[5]

\usepackage{amsmath}
\usepackage{amsfonts}
\newtheorem{thm}{Theorem}

\newdefinition{rmk}{Remark}
\newproof{pf}{Proof}
\newproof{pot}{Proof of Theorem \ref{thm2}}

\begin{document}

\begin{frontmatter}

\title{Satellite Conjunction Analysis and the False Confidence Theorem}

%% Group authors per affiliation:

\author{Michael Scott Balch\footnote{Alexandria Validation Consulting, LLC, 8565 Richmond Highway, \#203, Alexandria, VA 22309, United States}, Ryan Martin\footnote{Department of Statistics, North Carolina State University, 5109 SAS Hall, Raleigh, NC 27695, United States}, Scott Ferson,\footnote{Institute of Risk and Uncertainty, Liverpool University, Chadwick Building, Liverpool L69 7ZF, United Kingdom}}

%
%\author{Michael Scott Balch}
%\address{Alexandria Validation Consulting, LLC, 8565 Richmond Highway, \#203, Alexandria, VA 22309, United States}
%
%\author{Ryan Martin}
%\address{Department of Statistics, North Carolina State University, 5109 SAS Hall, Raleigh, NC 27695, United States}
%
%\author{Scott Ferson}
%\address{Institute of Risk and Uncertainty, Liverpool University, Chadwick Building, Liverpool L69 7ZF, United Kingdom}

\begin{abstract}

Satellite conjunction analysis is the assessment of collision risk during a close encounter between a satellite and another object in orbit. A counterintuitive phenomenon has emerged in the conjunction analysis literature, namely, probability dilution, in which lower quality data paradoxically appear to reduce the risk of collision. We show that probability dilution is a symptom of a fundamental deficiency in probabilistic representations of statistical inference, in which there are propositions that will consistently be assigned a high degree of belief, regardless of whether or not they are true. We call this deficiency false confidence. In satellite conjunction analysis, it results in a severe and persistent underestimate of collision risk exposure. 

We introduce the Martin--Liu validity criterion as a benchmark by which to identify statistical methods that are free from false confidence. Such inferences will necessarily be non-probabilistic. In satellite conjunction analysis, we show that uncertainty ellipsoids satisfy the validity criterion. Performing collision avoidance maneuvers based on ellipsoid overlap will ensure that collision risk is capped at the user-specified level. Further, this investigation into satellite conjunction analysis provides a template for recognizing and resolving false confidence issues as they occur in other problems of statistical inference.

\end{abstract}

\begin{keyword}
belief functions\sep epistemic uncertainty\sep statistical inference\sep collision risk\sep satellite conjunction analysis\sep probability dilution

% \texttt{elsarticle.cls}\sep \LaTeX\sep Elsevier \sep template
% \MSC[2010] 00-01\sep  99-00
\end{keyword}

\end{frontmatter}

%\linenumbers

\section{Introduction}
\label{sect:intro}

A satellite conjunction is an event in which two satellites or a satellite and a piece of debris are estimated to pass near each other. The goal of satellite conjunction analysis is to determine whether the risk of collision is high enough to necessitate a collision avoidance maneuver. Avoiding collision is not just important to the individual satellite operator. Everyone in the space industry has a stake in minimizing the number of in-orbit collisions. Each collision adds to the debris already in orbit, and in aggregate, this debris has the potential to make low Earth orbit unnavigable \cite{kessler2010kessler,liou2006risks,liou2008instability}. % The situation is further exacerbated by the fact that the number of active satellites is expected to greatly increase over the next few years \cite{India2017Launch,wapo2016floodorbit}. 
So, while some of the issues involved in conjunction analysis may initially seem arcane, correctly resolving those issues is of existential importance to the space industry.

\subsection{A Modern Problem}
\label{subsect:intro.modern}

Over the past 15 years, aerospace researchers have recognized a counterintuitive phenomenon in 
%the collision risk numbers obtained in 
satellite conjunction analysis, known as {\em probability dilution} \cite{alfano2003collision}. That is, as uncertainty in the satellite trajectories increases, the epistemic probability of collision eventually decreases. Since trajectory uncertainty is driven by errors in the tracking data, the seemingly absurd implication of probability dilution is that lower quality data reduce the risk of collision. 

Several researchers have attempted to synthesize an ad hoc fix to probability dilution \cite{alfano2003collision,frigm2009single,plakalovic2011tuned,sun2014spacecraft,balch2016corrector,carpenter2017relevance}. Most of them have defined an alternative collision risk metric that is increasing, or at least non-decreasing, as a function of trajectory uncertainty. %However, none of these researchers has offered a rigorous explanation for why probability dilution is supposedly wrong or, at least, how it is misleading. % They all proceed on the basis that its counterintuitiveness is proof of its inappropriateness.
Alternatively, some satellite operators have adopted a high sensitivity to small collision probabilities \cite{newman2014evolution}. However, if probability dilution is indicative of a fundamental problem with collision probabilities, it is not immediately obvious that being sensitive to small collision probabilities will resolve that problem.

\subsection{A Centuries-Old Argument}
\label{subsect:intro.oldArgument}

Satellite trajectory estimation is fundamentally a problem of statistical inference, and the confusion over probability dilution cuts to the heart of the long-standing Bayesian--frequentist debate in statistics; for a review, see \cite{Barnett_Comparative}. Satellite orbits are inferred using radar data, optical data, GPS data, etc.~which are subject to random measurement errors, i.e., noise. The resulting probability distributions used in conjunction analysis represent epistemic uncertainty, rather than aleatory variability. That is to say, it is the trajectory estimates that are subject to random variation, not the satellite trajectories themselves. In the Bayesian view, that is a distinction without a difference; it is considered natural and correct to assess the probability of an event of interest, such as a possible collision between two satellites, based on an epistemic probability distribution \cite{Laplace_Essay,Finetti_Logiques,robert2007bayesian}. In the frequentist view, however, it is considered an anathema to try to compute the probability of a non-random event \cite{Boole_Laws,neyman1937outline,Mayo_Book}. This prohibition is expressed in some corners of the uncertainty quantification community as the requirement that epistemic uncertainty be represented using non-probabilistic mathematics \cite{ferson1996different,Salicone_Text,Roy_Oberkampf_Book}. Despite these objections, at present, satellite navigators appear to be adopting the Bayesian position.

\subsection{Definitions}
\label{subsect:intro.defs}

% Since this paper is intended for a general scientific audience, it is necessary to define key terms from the uncertainty quantification literature. 

An \textit{aleatory uncertainty} is an uncertainty due to random (i.e., aleatory) variation of the quantity or event being analyzed. An \textit{epistemic uncertainty} is uncertainty due to imperfect knowledge about the quantity or event being analyzed. One of the many ways in which epistemic uncertainty can arise in a fixed quantity is if the data used to measure or infer it are subject to random errors.

A \textit{belief function} is a function used to represent the amount of positive support a proposition has accrued from the evidence. A \textit{plausibility function} merely represents the lack of evidence against a proposition. Belief and plausibility are complementary functions; that is, $\textnormal{Pls}\!\left( A \right) = 1 - \textnormal{Bel}\!\left( \neg A \right)$, for every proposition, $A$, and its complement (i.e., negation) $\neg A$. As a rule, the plausibility accorded to a proposition is always greater than or equal to the belief accorded to it. An additive belief function is a belief function that satisfies the Kolmogorov axioms and in which, as a consequence, $\textnormal{Bel}\!\left( A \right) = \textnormal{Pls}\!\left( A \right)$ for all propositions \cite{Shafer_1976_Book}. 

An \textit{epistemic probability} is an additive belief function used to represent epistemic uncertainty. Confidence distributions,\footnote{We refer to confidence distributions in the sense of \cite{Balch_Conf_2012}, as opposed to \cite{Schweder_Hjort_book} and \cite{Singh_2013_CD}, which eschew epistemic probabilities.} fiducial distributions, and most Bayesian posterior distributions\footnote{One could argue that a posterior probability derived using an aleatory prior, as in \cite{vonMises_1942_BayesFormula}, is itself an aleatory probability.} are epistemic probability distributions. An \textit{aleatory probability} is a probability used to represent aleatory uncertainty. It is equal to the fraction of trials in which a specified event occurs---or would occur---over an infinite number of equivalent random trials. 
The mathematics of aleatory and epistemic probability functions are identical; the distinction is a question of subject matter. %The use of aleatory probabilities is generally accepted; only epistemic probability is a topic of controversy. 

\subsection{Our Contribution}
\label{subsect:intro.contribution}

In this paper, we pursue basic questions raised by probability dilution. Can satellite operators take epistemic probability at face value as a risk metric? If not, why not? What alternatives are safe for use in conjunction analysis? How can we know that those alternatives are safe? We address these questions, in part, by placing conjunction analysis in its context as a problem of statistical inference. Consequently, what we discover about probability dilution has broad implications for statisticians as well as satellite operators. % These discoveries are expressed via three contributions that we make in Sections~\ref{sect:satellite}-\ref{sect:validity}. 

In Section~\ref{sect:satellite}, we provide a thorough exposition of 
%and explanation for 
probability dilution. 
%This has previously been lacking in the literature. 
We show that the counterintuitiveness of probability dilution is entirely due to the fact that, in conjunction analysis, probability of collision is an epistemic probability, rather than an aleatory probability. Moreover, we show that probability dilution is linked to a more severe operational defect that we call false confidence. We then evaluate the efficacy of pursuing small epistemic probabilities as a way to mitigate the problems posed by false confidence.

In Section~\ref{sect:falseConfidenceTheorem}, we show that false confidence is not limited to conjunction analysis. All epistemic probability distributions are subject to false confidence. That is, in any additive belief function used to represent epistemic uncertainty, there are false propositions that have a high probability of being assigned a high belief value. The false confidence generated by epistemic probability of collision in satellite conjunction analysis is a manifestation of this more general deficiency.  

In Section~\ref{sect:validity}, we demonstrate a framework for identifying statistical methods that are free from false confidence. We build this approach on the ``validity criterion'' introduced in the Martin--Liu Theory of Inferential Models \cite{martin2016book}. For a statistical inference to be valid in the Martin--Liu sense, only a true proposition may have a high probability of being assigned a high belief value. An inference satisfying this criterion will not suffer from the false confidence phenomenon described in Section~\ref{sect:falseConfidenceTheorem}. Note that the reader does not need to be familiar with Martin--Liu theory to understand the work presented in this paper. Only the validity criterion is used, and its formal definition is provided in Section~\ref{sect:validity}. % From there, we prove that confidence regions in general satisfy the Martin--Liu validity criterion. We then make the case that $K\sigma$ uncertainty ellipsoids around predicted satellite positions can usually be interpreted as confidence regions.

\section{Epistemic Probability of Collision}
\label{sect:satellite}

% Conjunction analysis consists of many small decisions and analyses culminating in one big decision: whether or not to perform a collision avoidance maneuver. From the time a satellite navigator is first made aware of an upcoming conjunction, there are many calculations of collision risk, updated as new data become available. Most satellite conjunctions are resolved by simply waiting for better data that conclusively show that there is not going to be a collision \cite{newman2016nasa}. 

% Nevertheless, collisions have happened \cite{newman2009notbigsky}. Moving forward, it is imperative that satellite navigators successfully detect and prevent additional collisions. That imperative will grow as the number of active satellites increases. A crucial element in the success or failure of the detection of an impending collision is the statistical reliability of the risk metrics that are used. 

Epistemic probability of collision is coming into more common use as a risk metric in satellite conjunction analysis \cite{newman2014evolution,alfano2016probability,frigm2015total}. In fact, the most-cited literature on satellite conjunction analysis simply takes it for granted that epistemic probability of collision is the correct way to represent collision risk, e.g., \cite{alfriend1999probability,patera2001general,chan2008spacecraft}. However, the counterintuitiveness of probability dilution calls this practice into question, especially considering the unsettled status of epistemic probability in the statistics and uncertainty quantification communities. 

This part of the paper explores issues surrounding the use of epistemic probability of collision as a risk metric for satellite conjunction analysis. Section~\ref{subsect:satellite.compute} reviews the standard computation of collision probability. Section~\ref{subsect:satellite.dilution} describes how probability dilution arises from those mathematics. Section~\ref{subsect:satellite.false} describes a more severe operational defect in epistemic collision probabilities: false confidence. Section~\ref{subsect:satellite.pursuing} explores the use of small thresholds for epistemic probability as a way to counteract false confidence. 

\subsection{Computation of Collision Probability}
\label{subsect:satellite.compute}

As mentioned in Section~\ref{subsect:intro.oldArgument}, satellite conjunction analysis is rooted in a problem of statistical inference. The inferred parameter of interest is the vector of positions and velocities of the two satellites at time of closest approach. That is, 
\[
\theta = \left( u_1 , v_1 , w_1 , \dot{u}_1 , \dot{v}_1 , \dot{w}_1 , u_2 , v_2 , w_2 , \dot{u}_2 , \dot{v}_2 , \dot{w}_2 \right)^\top
\]
where $u_1,v_1,w_1$ are the position of one satellite; $\dot{u}_1 , \dot{v}_1 , \dot{w}_1$ are its velocity; $u_2 , v_2 , w_2$ are the position of the other satellite; $\dot{u}_2 , \dot{v}_2 , \dot{w}_2$ are its velocity.  The true value of $\theta$ is unknown; it is estimated using tracking data, $x$, via an inferential algorithm called a ``filter'' that delivers an estimate, $\hat{\theta}( x )$, along with its uncertainty, expressed as a $12 \times 12$ covariance matrix, $C_{\Theta}$. In this paper, uncertainty in $\theta$ is attributed to random errors in the tracking data. In reality, uncertainty also arises from errors in the dynamics model used in the filter. 

Uncertainty in the trajectory estimates is usually assumed to have a multivariate normal distribution \cite{alfano2003collision,alfriend1999probability,patera2001general}. The epistemic probability density for $\theta$ is therefore taken to be
\[
f_{\Theta;x}(\theta) = \Bigl\{ \left( 2\pi \right)^{ 12 } \det( C_{\Theta} ) \Bigr\}^{-1/2} \exp{ \left\{ -\frac{1}{2} \left( \hat{\theta}( x ) - \theta \right)^\top C_{\Theta}^{-1} \left( \hat{\theta}( x ) - \theta \right) \right\} }
\]
where $\det( C_{\Theta} )$ is the determinant of $C_{\Theta}$. This epistemic probability distribution for $\theta$ can be rationalized in two ways. One way is to treat filter estimation as an implementation of Bayesian inference \cite{meinhold1983understanding}. The other way is to describe this representation as a fiducial distribution for $\theta$, in the sense of \cite{hannig2009generalized}, or a confidence distribution for $\theta$, in the sense of \cite{Balch_Conf_2012}, on the argument that a random resampled $\hat{\theta}\!\left( X \right)$ would have a multivariate normal distribution with mean $\theta$ and covariance matrix $C_{\Theta}$, where $X$ represents a random realization of the tracking data, drawn from the same distribution as $x$, the tracking data actually obtained. The opening of Section~\ref{subsect:validity.ellipsoids} briefly describes the conditions under which the normality assumption will and will not hold. 

Assuming normality holds, the standard epistemic probability of collision calculation proceeds in three steps. First, uncertainty in $\theta$ is propagated to uncertainty in the relative offset between the two satellites at closest approach. That is,
\[
\Delta u = u_2 - u_1 \ \ \ \ \Delta v = v_2 - v_1 \ \ \ \ \Delta w = w_2 - w_1.
\]
Since the linear transform of a multivariate normal variable is itself a multivariate normal, the triple $\left( \Delta u , \Delta v , \Delta w \right)$ has a multivariate normal distribution with the following $3\times3$ covariance matrix:
\[
C_{\Delta} = A C_{\Theta} A^\top \ \ \ \ \textnormal{where} \ \ \ \ A =
\begin{bmatrix}
-1 & 0 & 0 & 0 & 0 & 0 & +1 & 0 & 0 & 0 & 0 & 0 \\
0 & -1 & 0 & 0 & 0 & 0 & 0 & +1 & 0 & 0 & 0 & 0 \\
0 & 0 & -1 & 0 & 0 & 0 & 0 & 0 & +1 & 0 & 0 & 0
\end{bmatrix}
\]
or more compactly
\[
C_{\Delta} = C_{\Theta;1:3,1:3} + C_{\Theta;7:9,7:9} - 2 C_{\Theta;1:3,7:9},
\]
where $C_{\Theta;m_1:m_2,n_1:n_2}$ is the sub-matrix of all entries in rows $i$ and columns $j$ such that $m_1 \leq i \leq m_2$ and $n_1 \leq j \leq n_2$. That is, $C_{\Delta}$ is the sum of the covariance matrices for the positions of the two satellites minus the covariance between the positions of the two satellites.

In the second step, uncertainty along the axis parallel to the relative velocity vector is integrated out. That is, instead of computing probability of collision as the probability of the true value of the displacement falling within the set of $\left( \Delta u , \Delta v , \Delta w \right)$ values indicative of collision, probability of collision is computed as the probability accorded to an extruded ellipse containing that three-dimensional set. The theoretical rationale behind this move is that it captures uncertainty in the timing of closest approach \cite{patera2001general}. The primary practical effect of this step is to reduce a three-dimensional integration problem to a two-dimensional integration problem. More nuanced ways of treating uncertainty in the timing of closest approach are available \cite{hall2017time,helpSTK_nonLinearProb}, but the standard two-dimensional treatment is sufficient for our analysis. 

% The primary practical effect of this step is to reduce a three-dimensional integration problem to a two-dimensional integration problem in a way that obtains a probability of collision value that is always larger than the mathematically correct instantaneous probability of collision value. 

This transformation is executed in a few sub-steps. It starts with the unit vector in the relative velocity direction, $i _{\Delta V}$, defined as
\[
\Delta V = \left( \dot{u}_2 , \dot{v}_2 , \dot{w}_2 \right) - \left( \dot{u}_1 , \dot{v}_1 , \dot{w}_1 \right) \quad \textnormal{and} \quad i_{\Delta V} = \Delta V / \| \Delta V \|.
\]
Next, another direction in the plane perpendicular to $i_{\Delta V}$ can be defined arbitrarily. For example it could be taken as in the same direction as the cross-product of $i_{\Delta V}$ and $\left( \Delta u , \Delta v , \Delta w \right)$, assuming that they are not aligned. Denote this arbitrary direction vector as $i_{u'}$. The final direction vector is simply the cross-product, $i_{v'} = i_{\Delta V} \times i_{u'}$. Once this system has been established, define the rotated displacement vector as
\[
\begin{bmatrix} u' \\ v' \\ w' \end{bmatrix} = M \begin{bmatrix} \Delta u \\ \Delta v \\ \Delta w \end{bmatrix} \ \ \ \ \textnormal{where} \ \ \ \ M = \begin{bmatrix} i_{u'}^\top \\ i_{v'}^\top \\ i_{\Delta V}^\top \end{bmatrix}
\]
is the rotation matrix. The covariance matrix for the rotated $\left( u' , v' , w' \right)$ vector is $C_{\Delta'} = M C_{\Delta} {M}^{\top}$.  Since $\left( u' , v' , w' \right)$ has a multivariate normal distribution, uncertainty in the $\Delta V$ direction can be integrated out by simply dropping $w'$. That completes the extrusion step, reducing the problem to the two-dimensional $\left( u' , v' \right)$ space, with estimate $\left( \hat{u}' , \hat{v}' \right)$ and covariance matrix $C_{\Delta';1:2,1:2}$. 

The third and final step in the standard calculation of collision probability is to define the probability of collision as the epistemic probability of having $u'$ and $v'$ such that $u'^2 + v'^2 \leq R^2$, where $R$ is the sum of the characteristic radii of the two satellites in the conjunction event. In other words, the satellite shapes are 
%ignored; they are 
approximated as spherical objects. Collision occurs if and only if, at closest approach, the distance between the two satellites is less than their combined size.

This probability integral can be re-expressed in a standardized probability space, by performing an eigendecomposition of the displacement covariance matrix \cite{balch2016corrector}, which yields
\[
\begin{bmatrix} u' \\ v' \end{bmatrix} = \begin{bmatrix} \hat{u}' \\ \hat{v}' \end{bmatrix} - E_S \begin{bmatrix} S_1 & 0 \\ 0 & S_2 \end{bmatrix} \begin{bmatrix} \xi_1 \\ \xi_2 \end{bmatrix}
\]
where $S^2_1$ and $S^2_2$ are the eigenvalues of $C_{\Delta';1:2,1:2}$; $E_S$ is a $2 \times 2$ matrix whose columns are the normalized eigenvectors of $C_{\Delta';1:2,1:2}$; and finally, $\xi_1$ and $\xi_2$ are independent unit normal variables, with mean zero and variance one. Note that, since $E_S$ is the eigendecomposition of a symmetric matrix, $E^\top_S E_S$ returns the identity matrix. Defining $u''$ and $v''$ as the $E^\top_S$ rotation of $u'$ and $v'$, the magnitude of $\left( u'' , v'' \right)$ is equal to the magnitude of $\left( u' , v' \right)$. That is,
\[
\begin{bmatrix} u'' \\ v'' \end{bmatrix} = E^\top_S \begin{bmatrix} u' \\ v' \end{bmatrix} \implies  u''^2 + v''^2 = u'^2 + v'^2, \quad \textnormal{where} \quad \begin{cases} u'' = \hat{u}'' - S_1 \xi_1 \\ u'' = \hat{v}'' - S_2 \xi_2. \end{cases} 
%u'' = \hat{u}'' - S_1 \xi_1, \quad \textnormal{and} \quad v'' = \hat{v}'' - S_2 \xi_2.
\]
Thus, the following equivalences hold for the integration condition:
\[
u'^2 + v'^2 \leq R^2 \iff u''^2 + v''^2 \leq R^2 \iff \left( S_1 \xi_1 - \hat{u}'' \right)^2 + \left( S_2 \xi_2 - \hat{v}'' \right)^2 \leq R^2.
\]
What has been accomplished is that now, instead of integrating a circle in a $\left( u' , v' \right)$ space that will vary from problem to problem, we are integrating an ellipse in a standardized unit normal $\left( \xi_1 , \xi_2 \right)$ space. The centroid of that ellipse is located at $\left( \hat{u}'' / S_1 , \hat{v}'' / S_2 \right)$, and its semi-major and semi-minor axes are $R/S_1$ and $R/S_2$.

Following the logic of \cite{patera2005calculating}, the epistemic probability of collision is computed using a contour integral transformation, as follows:
\begin{align*}
\textnormal{Bel}\!\left( \mathbb{C} \right) & = \iint_{ \left( S_1 \xi_1 - \hat{u}'' \right)^2 + \left( S_2 \xi_2 - \hat{v}'' \right)^2 \leq R^2 } { \frac{1}{2\pi} \exp\!{ \left[ -\frac{1}{2}\left( \xi^2_1 + \xi^2_2 \right) \right] } }d\xi_1 d\xi_2 \\
& = \int_0^{2\pi}{ \left( \frac{1-e^{-r_\psi^2/2}}{2 \pi r_\psi^2} \right) \left( \frac{R^2}{S_1 S_2} + \frac{\hat{u}'' R}{S_1 S_2} \cos\!\psi + \frac{\hat{v}'' R}{S_1 S_2} \sin\!\psi \right)  }d\psi
\end{align*}
where 
\[
r_\psi^{2} = \left( \frac{\hat{u}''}{S_1} + \frac{R}{S_1} \cos\!\psi \right)^2 + \left( \frac{\hat{v}''}{S_2} + \frac{R}{S_2} \sin\!\psi \right)^2.
\]
This integral can be approximated to negligible error using the trapezoidal rule, provided that at least $10 \times \max\left( S_1 / S_2 , S_2 / S_1 \right)$
evenly spaced quadrature points are used \cite{balch2016corrector}.

Because our goal is to understand and explore how the epistemic probability of collision changes for different levels of uncertainty, we pursue the special case where $S = S_1 = S_2$. Defining the estimated displacement as $D^2 = \hat{u}''^2 + \hat{v}''^2 = \hat{u}'^2 + \hat{v}'^2$, this simplification reduces the epistemic probability of collision to a function of two ratios: estimated relative displacement, $D / R$, and relative uncertainty, $S / R$.

\subsection{Probability Dilution}
\label{subsect:satellite.dilution}

\begin{figure}[t]
	\centering
	\includegraphics[width=0.49\textwidth]{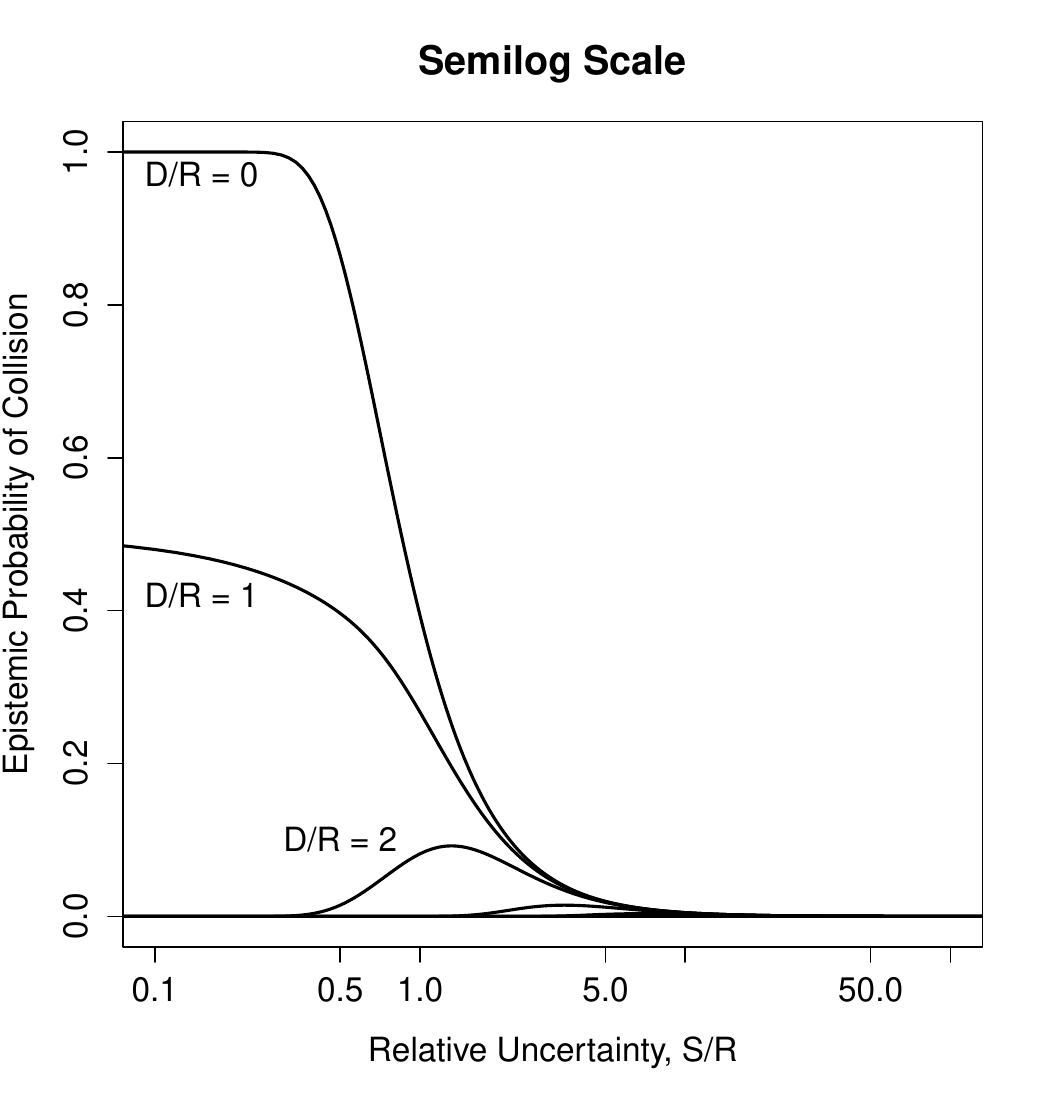}
	\includegraphics[width=0.49\textwidth]{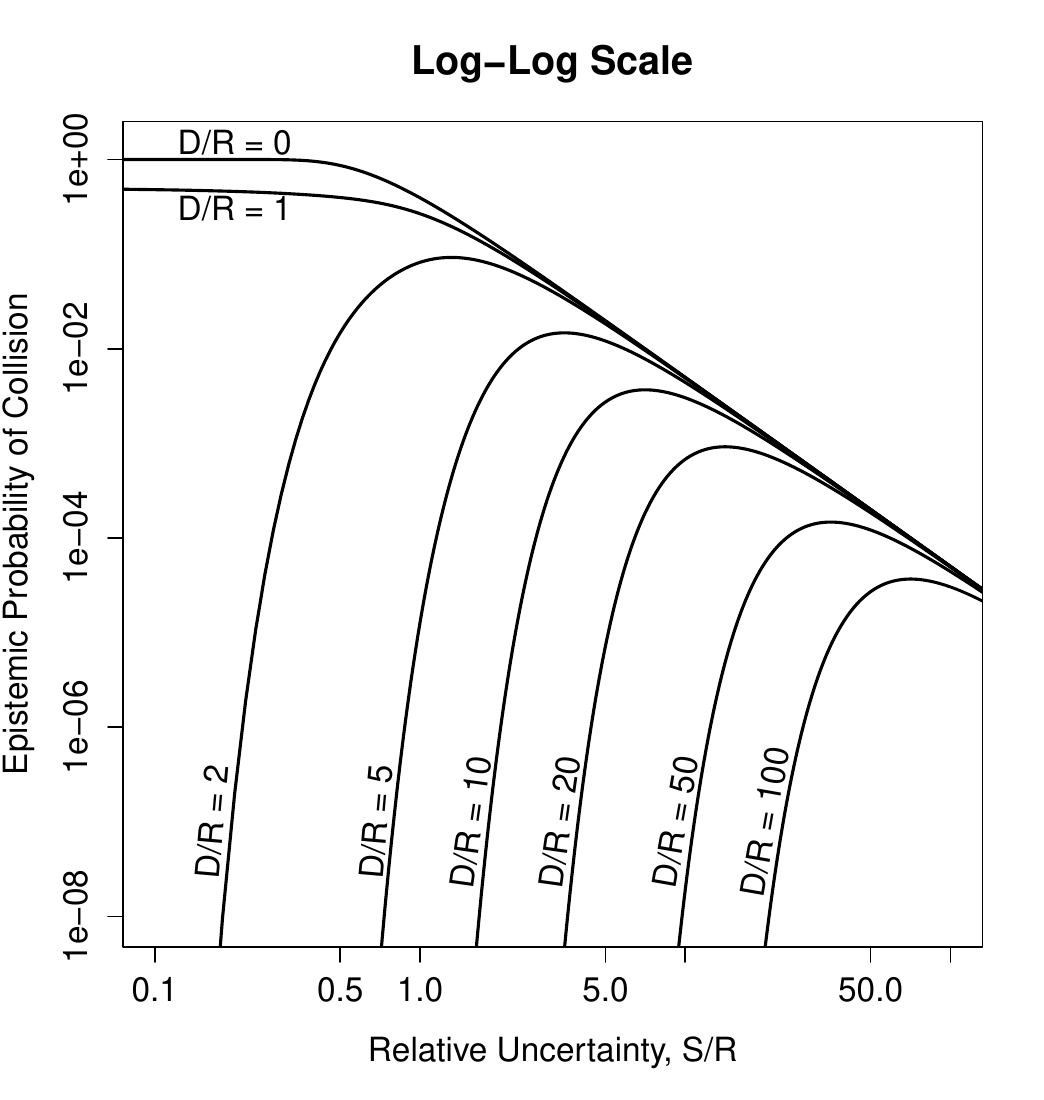}
	\caption{Epistemic probability of collision as a function of estimated relative displacement, $D / R$, and relative uncertainty, $S / R$, for the special case $S = S_1 = S_2$.}
	\label{fig:CollisionProbability}
\end{figure}

Figure~\ref{fig:CollisionProbability} illustrates epistemic probability of collision as a function of $S / R$ for several values of $D / R$. These curves follow a common pattern. So long as $D/R > 1$, for small uncertainties, the probability of collision is small. That much makes intuitive sense. If the satellites are estimated to miss each other, and the uncertainty in those trajectory estimates is small, then it is natural that the analyst have a high confidence that the satellites are not going to collide. Similarly, it makes sense that as uncertainty grows, the risk of collision also grows. However, something odd happens in the curves in Figure~\ref{fig:CollisionProbability}. Epistemic probability of collision eventually hits a maximum, and past that maximum, as relative uncertainty rises, the epistemic probability of collision {\em decreases}. 

This decrease is called probability dilution, and it has an odd implication. Since the uncertainty in the estimates of the trajectories reflects the limits of data quality, probability dilution seems to imply that lowering data quality makes satellites safer. That implication is counterintuitive in the extreme \cite{alfano2003collision,frigm2009single,plakalovic2011tuned,sun2014spacecraft,balch2016corrector}. As a rule, lowering the data quality makes any engineering system less safe, and to claim that ignorance somehow reduces collision risk seems foolish on its face. 

\begin{figure}[t]
	\centering
	\includegraphics[width=0.49\textwidth]{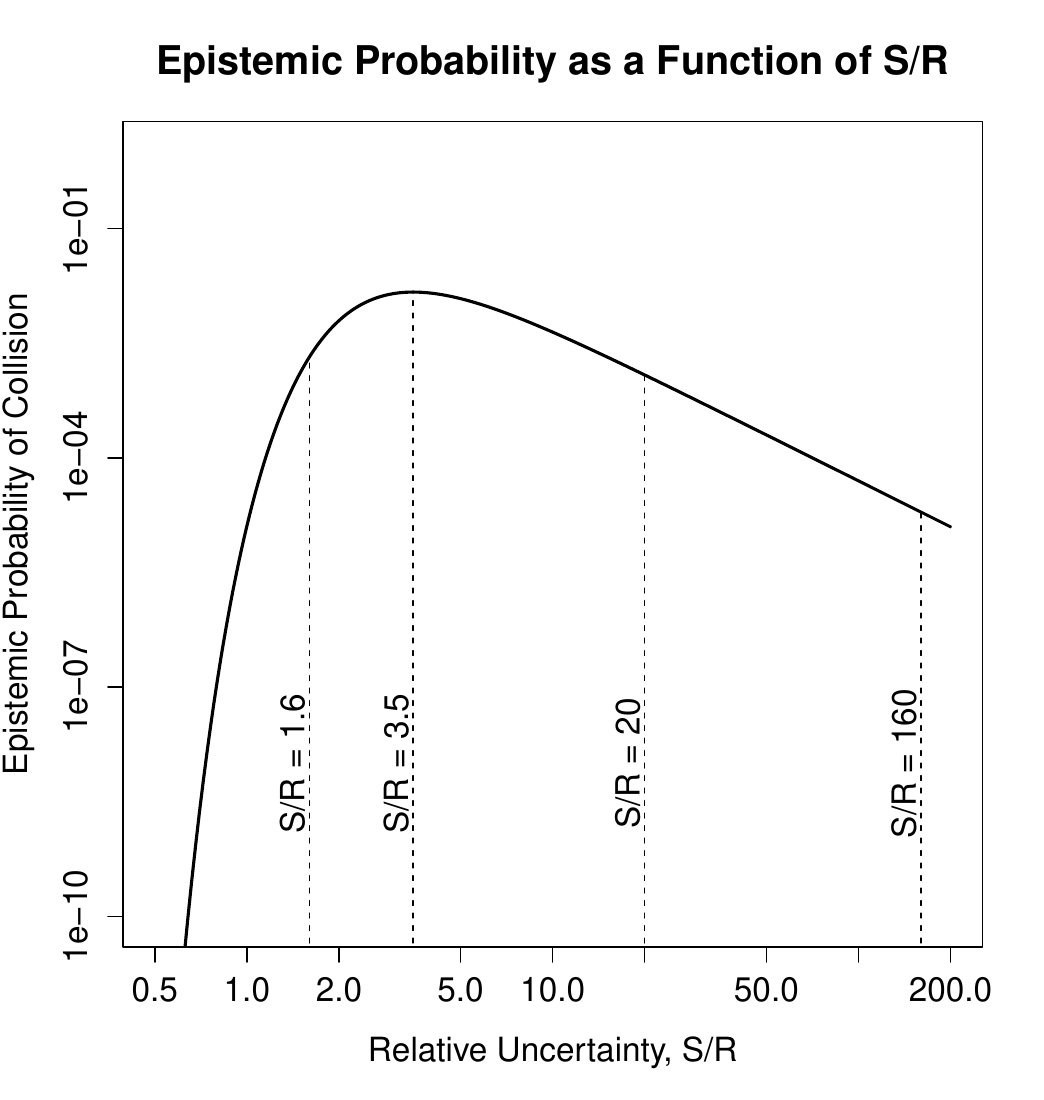}
	\includegraphics[width=0.49\textwidth]{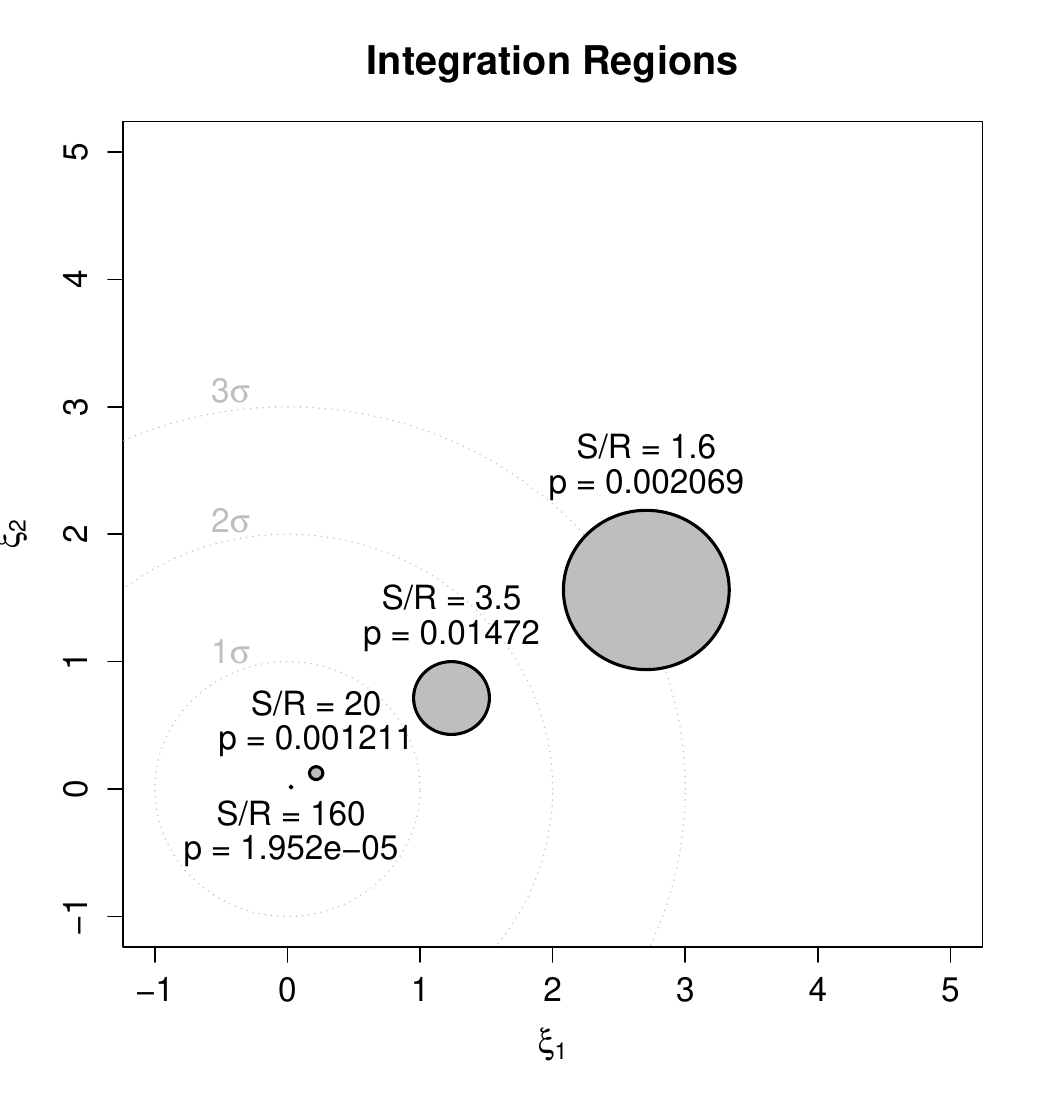}
	\caption{Epistemic probability of collision and integration region, given $D/R = 5$ and $S/R \in \left\{ 1.6 , 3.5 , 20 , 160 \right\}$. In the right panel, the $1\sigma$, $2\sigma$, and $3\sigma$ contours for unit normal space are marked with dashed lines for reference. These three contours capture $39.3\%$, $86.5\%$, and $98.9\%$ of the two-dimensional unit normal probability mass, respectively.}
	\label{fig:ShrinkingProbability}
\end{figure}

That being said, the mathematics of probability dilution are ironclad and relatively simple. Recall that, when expressed in unit normal space, the integration region determining the probability of collision is an ellipse whose centroid is at $\left( \hat{u}'' / S_1 , \hat{v}'' / S_2 \right)$ and whose semi-major and semi-minor axes are $R/S_1$ and $R/S_2$. Figure~\ref{fig:ShrinkingProbability} illustrates how the integration region changes as a function of $S$, when $S = S_1 = S_2$. Combined satellite size, $R$, and estimated distance at closest approach, $D$, are held fixed with $D/R = 5$. For four different values of $S/R$, the integration region is illustrated, along with the corresponding probability of collision. 

Increasing the level of uncertainty has two effects on the integration region: It shrinks and draws closer to the origin.  
%The integration region shrinks, and it draws closer to the origin. 
The probability density of a unit normal space is highest at the origin, decreasing exponentially as a function of the square of distance from the origin. If $D / R > 1$, the positive effect due to getting closer to the origin, where there is more probability mass, initially outweighs the negative effect due to shrinking the integration region. Near the origin, though, the probability density curve flattens out. So, after a certain point, there is not much more probability mass to be gained by moving closer to the origin, and the negative effect of shrinking the integration region overtakes the positive effect due to shifting the location of the integration region. Past that point, probability of collision decreases as uncertainty increases. So, in total, probability dilution is 
%That is why probability dilution happens, from a purely mathematical perspective; it is 
due to a simple and straightforward shrinkage of the integration region relative to displacement uncertainty.

Even if the normality and shape assumptions were relaxed, the integration region would still undergo the same shrinkage phenomenon that is illustrated in Figure~\ref{fig:ShrinkingProbability}. Each satellite represents a bounded set of points. The displacement values indicative of collision are therefore also a bounded set of points. If one were to relax the assumptions outlined in Section~\ref{subsect:satellite.compute}, instead of an ellipse, one would have a different shrinking shape in a standardized three-dimensional probability space. Probability of collision would still have the same qualitative behavior as a function of relative uncertainty. Mathematically speaking, the key factor underpinning probability dilution is that the integration region (i.e.,~failure domain) is bounded.

There is no way of reframing satellite conjunction analysis so that the failure domain is not bounded. Therefore, if we take the view that epistemic probability is valid, then it would seem that probability dilution is fundamental to conjunction analysis. Why, then, do aerospace researchers \cite{alfano2003collision,frigm2009single,plakalovic2011tuned,sun2014spacecraft,balch2016corrector} find this supposed mathematical inevitability so counterintuitive?

Interestingly, were the uncertainty in the satellite trajectories aleatory, rather than epistemic, probability dilution would actually make sense. For example, suppose two satellites were known, with certainty, to be on a collision course, and the satellite operator could only impart an impulse of random magnitude in a random direction, say, via a poorly controlled maneuvering thruster on a tumbling satellite. If the mean of this distribution were the null vector, due to a lack of control over thrust direction, then the higher the variance of the added impulse, the bigger the resulting perturbation and hence the smaller the resulting probability of collision. In this example, the mean of the resulting trajectory distribution would still have the satellites on a collision trajectory, but if one applies a big enough impulse to one of the satellites in a random direction, that poorly controlled collision avoidance maneuver still has a high probability of success. In that context, higher variance in a satellite's trajectory really does reduce the risk of collision. 

However, in this hypothetical example, it is a variance in the trajectory itself that makes the satellite safer, not a variance in the estimate of that trajectory. It should go without saying that, given two satellites on a sure collision trajectory, simply recomputing the trajectories with lower quality data does not make them safer. There is, therefore, an operational distinction to be made between genuine aleatory variability and epistemic uncertainty. 

To summarize, the problem with probability dilution is not the mathematics. The mathematics are incontrovertible, and they even make sense in the right context. Probability dilution is only counterintuitive when orbit ``variance'' reflects the limits of data quality, rather than genuine aleatory variability in the orbits themselves. So, if probability dilution is inappropriate, that inappropriateness must be rooted in a mismatch between the mathematics of probability theory and the epistemic uncertainty to which they are applied in conjunction analysis.

\subsection{False Confidence}
\label{subsect:satellite.false}

The final clue indicating the inappropriateness of probability dilution is the fact that, for a fixed $S/R$ ratio, there is a maximum computable epistemic probability of collision. Whether or not the two satellites are on a collision course, no matter what the data indicate, the analyst will have a minimum confidence that the two satellites will not collide. That minimum confidence is determined purely by the data quality.

% This is, even if the data indicates the satellites are on a direct collision course, there is a non-trivial upper bound on the probability of collision, which can be very small when the data is very noisy.

% % ; it is obtained even when the two satellites really are going to collide. 

For all $S/R$, the maximum collision probability is obtained when $D/R = 0$; that is, when the best estimate of the two satellite paths indicates that they are on a collision course. Therefore, the curve for $D/R = 0$ in Figure~\ref{fig:CollisionProbability} yields the maximum computable probability of collision as a function of $S/R$. For example, if the uncertainty in the distance between two satellites at closest approach is ten times the combined size of the two satellites, the analyst will {\em always} compute at least a $99.5\%$ confidence that the satellites are safe, even if, in reality, they are not. 

A relative uncertainty of ten may sound high, but it may help the reader to keep in mind the magnitude of the trajectories being estimated, relative to the size of the objects on those trajectories. Satellite trajectories are measured in thousands of kilometers; satellites are measured in meters. Trajectory uncertainties in conjunction analysis are usually measured in hundreds of meters \cite{sabol2010linearized}. For conjunction analysis done more than a week in advance, those uncertainties can grow to kilometers \cite{ghrist2012impact}. As a consequence, $S/R$ ratios greater than ten are the rule, not the exception. Under those circumstances, even if two satellites are on a collision course, because of probability dilution, the epistemic probability of collision is guaranteed to be small. 

The biggest problem with probability dilution is not its counterintuitiveness; instead, it is the \emph{false confidence} that it imbues. Due to probability dilution, satellite navigators using epistemic probability of collision as their risk metric will always compute a high confidence that their satellites are safe, regardless of whether or not they really are. The only rare exceptions are conjunctions in which both satellites are extremely well-tracked and at least one of them is relatively large (for an example, see \cite{newman2010overview}). From an operational perspective, if taken at face value, epistemic probability of collision is a misleading risk metric.

\subsection{Pursuing Small Epistemic Probabilities}
\label{subsect:satellite.pursuing}

%As mentioned in Section~\ref{subsect:intro.modern}, there are 
Some satellite navigators treat relatively small epistemic probabilities of collision as though they reflect a high level of risk. So, it would be inaccurate to characterize the entire satellite operations community as taking epistemic probabilities at face value. For example, the CARA group based at NASA Goddard monitors all conjunctions with a computed probability of collision greater than $10^{-7}$, and conjunctions with a computed probability of collision greater than $4.4 \times 10^{-4}$ are treated as high risk. These thresholds are calibrated in an effort to achieve desired frequentist error rates, in essence treating epistemic probability of collision as a test statistic, using historical data from past conjunctions as a reference population \cite{newman2014evolution}.

Unfortunately, even the idea of compensating for false confidence by setting a low threshold for epistemic probability of collision is problematic, because relative data quality modulates the relationship between epistemic probability of collision and the real aleatory probability of failing to detect an impending collision. To model this relationship, we consider the aleatory variability of predicted distance, $D$, between two satellites at closest approach, if the data used to produce that prediction were redrawn from the same distribution as the actual data obtained. Assuming that variability in the estimated direction of relative velocity is small enough to be neglected, the marginalization and eigendecomposition steps can proceed as in Section~\ref{subsect:satellite.compute}. This yields 
\[
D^2 = \hat{U}''^2 + \hat{V}''^2 =  \left( u'' + S_1 \xi_1 \right)^2 + \left( v'' + S_2 \xi_2 \right)^2
\]
where $\left( u'' , v'' \right)$ is the true unknown two-dimensional relative displacement between the two satellites; $\hat{U}''$ and $\hat{V}''$ are the random resampled analog to $\hat{u}''$ and $\hat{v}''$ from Section~\ref{subsect:satellite.compute}; $S_1$ and $S_2$ are the principal standard deviations of the marginalized two-dimensional covariance matrix, as in Section~\ref{subsect:satellite.compute}; and $\xi_1$ and $\xi_2$ are independent unit normal random variables. 

In the case where $S_1 = S_2 = S$, this simplifies to
\[
\left( \frac{D}{S} \right)^2 = \left( \frac{u''}{S} + \xi_1 \right)^2 + \left( \frac{v''}{S} + \xi_2 \right)^2.
\]
This can be further reduced to $(D / S)^2 = \chi^2_{2;\left( D_T/S \right)^2}$, where $\chi^2_{j,\delta}$ is a non-central chi-squared random variable with $j$ degrees of freedom and non-centrality parameter $\delta$, and $D_T^2 = u''^2 + v''^2$ is the unknown true distance between the two satellites at closest approach. 

In Section~\ref{subsect:satellite.compute}, we showed how epistemic probability of collision can be computed as a function of $D/S$ and $S/R$. Given a fixed $S/R$, epistemic probability of collision is a decreasing function of $D/S$. So, given a fixed detection threshold for epistemic probability, the aleatory probability of detecting an impending collision is equal to the probability of obtaining an observed $D/S$ that is less than or equal to the critical $D/S$ that yields an epistemic probability of collision equal to the threshold value. Rather than try to compute this relationship directly, it is easier to simply generate $D/S$ values from its non-central chi-squared distribution and then compute the corresponding epistemic probabilities of collision. This yields epistemic probability threshold as a function of desired detection rate, which can be reversed to give the aleatory probability of failing to detect an impending collision as a function of a satellite operator's chosen detection threshold for epistemic probability of collision. 

\begin{figure}[t]
	\centering
	\includegraphics[width=0.49\textwidth]{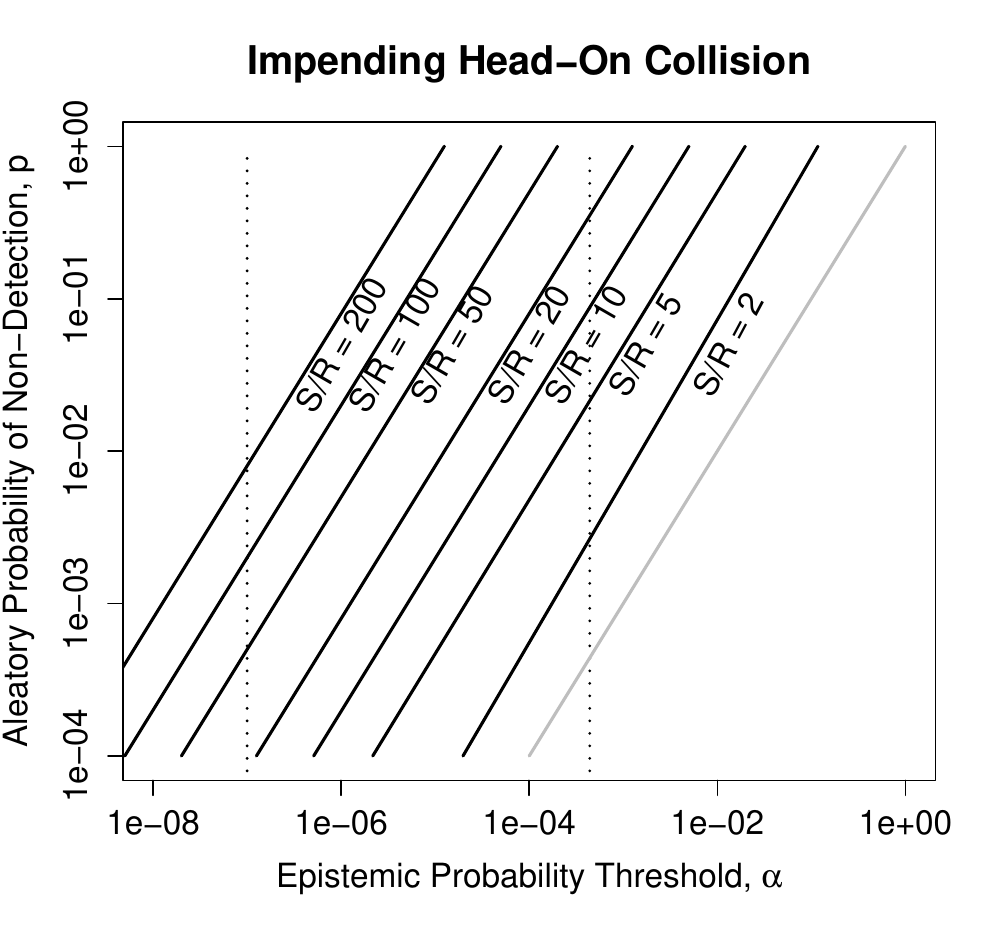}
	\includegraphics[width=0.49\textwidth]{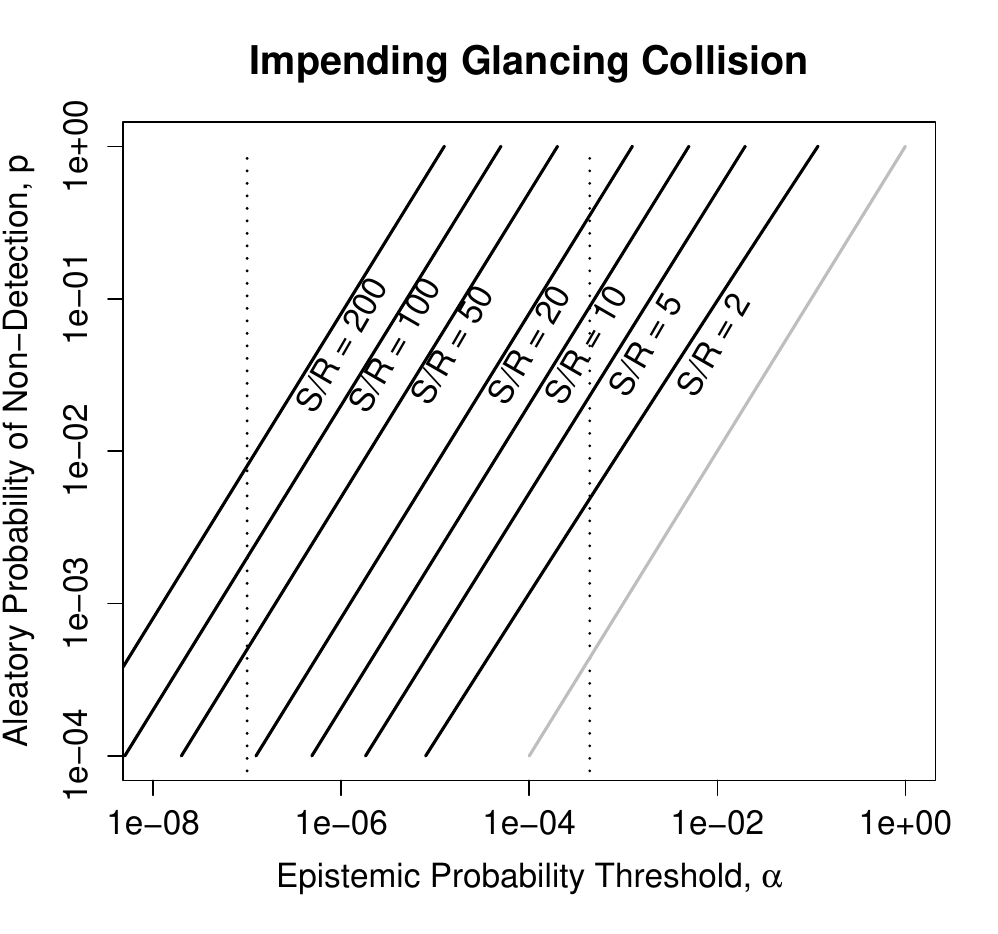}
	\caption{Aleatory probability of failing to detect an impending collision as a function of epistemic probability threshold for $S/R \in \left\{ 2 , 5 , 10 , 20 , 50 , 100 , 200 \right\}$. The left panel gives the error rate curves for a head-on collision, $D_T = 0$. The right panel gives the error rate curves for a glancing collision, $D_T = R$. Equivalently, these plots yield the aleatory probability, $p$, of obtaining data such that non-collision will be assigned an epistemic probability of at least $1-\alpha$, given the underlying reality that a collision is actually going to occur. The vertical dotted lines correspond to the CARA lower and upper thresholds reported in \cite{newman2014evolution} as $10^{-7}$ and $4.4 \times 10^{-4}$. The gray line in each panel demarcates equality between the nominal epistemic probability threshold and the real aleatory probability of failing to detect an impending collision.}
	\label{fig:CARA-Pc}
\end{figure}

Figure~\ref{fig:CARA-Pc} illustrates this relationship for several values of $S/R$ in two scenarios, an impending head-on collision with $D_T = 0$ and an impending glancing collision with $D_T = R$. The two dotted vertical lines mark the lower and upper thresholds for epistemic probability of collision reported in \cite{newman2014evolution}. Conjunctions with an epistemic probability of collision less than $10^{-7}$ are treated as lowest risk and ignored. Conjunctions with an epistemic probability of collision greater than $4.4 \times 10^{-4}$ are treated as high risk events. Using these thresholds, if $S/R = 10$, an impending head-on collision has a $91.2\%$ chance of correctly being identified as a high risk conjunction. Similarly, a glancing collision with $S/R = 10$ has a $91.1\%$ chance of being correctly identified as high risk. So, broadly speaking, one could say that the CARA threshold has reasonably good performance for $S/R = 10$. Unfortunately, that performance degrades quickly as uncertainty increases. At $S/R = 20$, the aleatory probability of correctly identifying an impending collision as a high risk conjunction drops to $64.8\%$. For $S/R > 33.74$, the aleatory probability of detecting an impending collision using the $4.4 \times 10^{-4}$ threshold drops to zero. For example, in a conjunction between two objects with a combined radius of five meters, a satellite navigator using the CARA upper risk threshold will have zero chance of detecting an impending collision unless or until their predicted displacement uncertainty drops below $170$ meters.

In summary, trying to compensate for false confidence by pursuing small epistemic probabilities of collision appears to be an unproductive strategy. Even small epistemic probabilities end up distorted, and they still result in false confidence under conditions that satellite navigators regularly encounter. Furthermore, the relationship between data quality and false confidence is more complicated than we have depicted in this idealized analysis. It depends not only on $S/R$ but also on $S_1/S_2$, and in reality, $S_1$ and $S_2$ are never equal. Every conjunction analysis has a unique combination of covariance matrices. So, if an analyst wants to hold the risk of failed detection at a consistent level, it is necessary to define a different epistemic probability threshold for each and every conjunction analysis. In other words, epistemic probability of collision could, theoretically, be used as a test statistic,
%\footnote{Setting $S = \sqrt{S_1 S_2}$, our analysis may provide a reasonable back-of-the-envelope estimate for the relationship between epistemic probability threshold and failed detection rate, so long as $S_1$ and $S_2$ are not too unequal.} 
but it cannot be taken directly as a risk metric, even if satellite navigators are sensitive to small epistemic probabilities of collision. There is no single threshold for epistemic probability that will detect impending collisions with a consistent degree of statistical reliability.

\section{The False Confidence Theorem}
\label{sect:falseConfidenceTheorem}

What causes probability dilution and false confidence? As explored in Section~\ref{subsect:satellite.dilution}, it is not an error in the mathematics. Rather, it appears to be a mismatch between the mathematics of probability theory and the subject matter to which those mathematics are applied in satellite conjunction analysis. As mentioned in Section~\ref{subsect:intro.oldArgument}, questions over the appropriateness of epistemic probability have never been settled. As shown in Section~\ref{subsect:satellite.pursuing}, those questions take on an undeniably practical dimension in satellite conjunction analysis. 

This part of the paper introduces the false confidence theorem. Section~\ref{subsect:fct.insight} suggests that probability dilution and false confidence are rooted in the axioms of probability theory. Section~\ref{subsect:fct.definition} posits a formal definition for false confidence. Section~\ref{subsect:fct.theorem} proves that all epistemic probability distributions assign copious amounts of false confidence. Section~\ref{subsect:fct.implications} places this theorem in the context of the Bayesian--frequentist debate.

\subsection{Motivating Insight}
\label{subsect:fct.insight}

Given highly uncertain data, it makes intuitive sense that one would not be certain of collision. An analyst should be reluctant to assign a high degree of belief to any proposition, based on low-quality data. The problem is that, within the confines of probability theory, general non-commitment of belief is not an option. One of the axiomatic properties of probability functions is that they are additive. So, any belief not assigned to ``collision'' is automatically assigned to ``not-collision,'' i.e., $\textnormal{Bel}\!\left( \neg \mathbb{C} \right) = 1 - \textnormal{Bel}\!\left( \mathbb{C} \right)$.  While a low assignment of belief to ``collision'' might seem natural when uncertainty is high, the consequently high assignment of belief to ``not-collision'' seems unnatural. The implication here is profound. If additivity is the root cause of false confidence, it means that false confidence is inherent to epistemic probability distributions, due to their mathematical structure. In Section~\ref{subsect:fct.theorem}, we show that this is exactly the case. 

\subsection{Definition of False Confidence}
\label{subsect:fct.definition}

Strictly speaking, one could characterize any confidence or belief or epistemic probability assigned to a false proposition as ``false confidence.'' However, the whole point of statistical inference is that one does not know, at the outset, the exact true value of the parameter being inferred. So, it is to be expected that, in the course of any given inference, some false propositions will be assigned some amount of belief. The problem described in Sections~\ref{subsect:satellite.false}-\ref{subsect:satellite.pursuing} is more severe: There is a fixed proposition of practical interest that is guaranteed or nearly guaranteed to be assigned a large epistemic probability, regardless of whether or not it is true. %high belief value.

This suggests a formal definition. To state this definition precisely, we need some additional notation.  First, let $\textnormal{Pro}_{X;\theta}$ denote the aleatory distribution of the data $X \in \Re^n$ depending on the fixed but unknown true value of the parameter $\theta \in \Omega_{\theta} \subseteq \Re^m$ being inferred.  Consider a set-function $\textnormal{Bel}_{\Theta;x}(\cdot)$ defined on the power set of $\Omega_{\theta}$, depending on data $x$, such that $\textnormal{Bel}_{\Theta;x}(A)$ represents the data analyst's epistemic degree of belief in the truthfulness of the proposition $A \subseteq \Omega_{\theta}$ about the unknown parameter $\theta$.  Then the map $x \mapsto \textnormal{Bel}_{\Theta;x}$ is an inferential method. We say that the inference about some proposition of interest, $A$, suffers from severe false confidence if, for some unacceptably high belief assignment $1-\alpha$ where $\alpha \in (0,1)$ and unacceptably high probability of assignment $p \in (0,1)$, there is some putative value of $\theta \in \Omega_{\theta}$ such that
\begin{equation}
\label{eq:fc}
A \not\ni \theta \quad \text{and} \quad \textnormal{Pro}_{X;\theta}( \{ x : \textnormal{Bel}_{\Theta;x}\!( A ) \geq 1 - \alpha \} ) \geq p. 
\end{equation}
Under this notation, Figure~\ref{fig:CARA-Pc} is a straightforward plot of $p$ as a function of $\alpha$, where $A$ is the assertion that the two satellites will not collide. As relative uncertainty grows, the level of false confidence seen in satellite conjunction analysis grows arbitrarily severe. Next, we show that, when considered over the full range of measurable propositions that the analyst might assess, most epistemic probability distributions used in practice suffer from arbitrarily severe false confidence. 

% and we say that it suffers from false confidence if, for any $\theta \in \Theta$, any $\alpha \in (0,1)$, and any $p \in (0,1)$, there exists a proposition $A \subset \Theta$ such that 

\subsection{Theorem}
\label{subsect:fct.theorem}

Most statistical inference in science and engineering involves one or more real-valued continuous parameter(s). 
%In general, let the observable be $X \in \Re^{n}$, and let the parameter be $\theta \in \Omega_\theta \subseteq \Re^{m}$ where $m$ is the number of parameters being inferred. 
So the epistemic probability distributions (e.g., posteriors) commonly used in practice are, for all values of the observable, continuous distributions over the parameter being inferred. That is, belief assigned by the additive belief function $\text{Bel}_{\Theta;x}$ can be represented via an epistemic probability density function, say, $f_{\Theta;x}(\theta)$, with respect to Lebesgue measure $\lambda$ on $\Omega_\theta$, depending on the observation, $x$. This is the case for most Bayesian posteriors, confidence distributions \cite{Balch_Conf_2012,Singh_2013_CD,Schweder_Hjort_book}, and fiducial distributions \cite{hannig2016generalized,Fisher_1935_Fiducial_Argument} used in practice. The satellite collision problem described in Section~\ref{sect:satellite} is one such example.  
%For example, in Section~\ref{subsect:satellite.compute}, epistemic probability of collision is derived from a multivariate normal distribution on the satellite states. 
The {\em false confidence theorem} below states that all such epistemic probability distributions suffer from arbitrarily severe false confidence. 

\begin{thm}
\label{thm:fc}
Consider an additive belief function $\textnormal{Bel}_{\Theta;x}$ characterized by an epistemic probability density function $f_{\Theta;x}$ on $\Omega_\theta $ such that $\sup_\vartheta f_{\Theta;x}(\vartheta) < \infty$ for $\textnormal{Pro}_{X;\theta}$-almost all $x$, for all $\theta$.  Then for any $\theta \in \Omega_\theta$, any $\alpha \in (0,1)$, and any $p \in (0,1)$, there exists a set $A \subset \Omega_\theta$ such that \eqref{eq:fc} holds.
\end{thm}

\begin{pf}
%The core strategy of this proof is to define a neighborhood, $B$, around the true parameter value and to then define $A$ as the complement of $B$. If $B$ is sufficiently small, then $A$ will be accorded the specified confidence at the specified probability, despite not containing the true parameter value.  
For fixed $x$, set $\bar f_x = \sup_\vartheta f_{\Theta;x}\!(\vartheta)$.  Then for any bounded set $B \subset \Omega_\theta$, 
\begin{equation}
\label{eq:bound}
\textnormal{Bel}_{\Theta;x}(B) = \int_B f_{\Theta;x}\!\left( \vartheta \right) \, d\vartheta \leq \bar f_x \, \lambda(B).
\end{equation}
As a function of $X \sim \textnormal{Pro}_{X;\theta}$, for fixed $\theta$, $\bar f_X$ has a distribution; let $\eta=\eta_{\theta,p} \in (0,\infty)$ 
%denote the $p^{\text{th}}$ quantile of the distribution of $\bar c_X$, i.e., 
satisfy 
\[ \textnormal{Pro}_{X;\theta}(\{x: \bar f_x \leq \eta\}) \geq p. \]
Choose a neighborhood $B \ni \theta$ with measure $\lambda(B) = \alpha / \eta$, so that $\lambda(B) \, \eta = \alpha$.  Define $A$ as the complement of $B$. For any belief function satisfying the Kolmogorov axioms, we have
\begin{equation}
\label{eq:3rdKolmo}
\textnormal{Bel}_{\Theta;x}(A) = 1-\textnormal{Bel}_{\Theta;x}(B).
\end{equation}
By definition, $A \not\ni \theta$, and by \eqref{eq:bound} and \eqref{eq:3rdKolmo}, we have 
\[ 
\bar f_x \leq \eta \implies \textnormal{Bel}_{\Theta;x}(B) \leq \lambda(B) \, \eta \iff \textnormal{Bel}_{\Theta;x}(B) \leq \alpha \iff \textnormal{Bel}_{\Theta;x}(A) \geq 1-\alpha. 
\]
By definition of $\eta$, the left-most event occurs with $\textnormal{Pro}_{X;\theta}$-probability at least $p$ and, therefore, the right-most event occurs with $\textnormal{Pro}_{X;\theta}$-probability at least $p$, proving the claim. 
\end{pf}

Theorem~\ref{thm:fc} is an existence result; so, our proof proceeds by constructing the simplest possible example. This is achieved by defining a neighborhood around the true parameter value that is so small that its complement---which, by definition, represents a false proposition---is all but guaranteed to be assigned a high belief value, simply by virtue of its size. In practice, no one would intentionally seek out such a proposition, but that is beside the point. 

Every real-world risk analysis problem involves a proposition of interest that is determined by the structure of the problem itself; e.g., ``Will these two satellites collide?''. Just as the practitioner will not seek out propositions strongly affected by false confidence, neither do practitioners have the option of avoiding such propositions when they arise. What the false confidence theorem shows is that, in most practical inference problems, there is no theoretical limit on how severely false confidence will manifest itself in an epistemic probability distribution, or more precisely, there is no such limit that holds for all measurable propositions. Such a limit can only be found for a specific proposition of interest through an interrogation of the belief assignments that will be made to it over repeated draws of the data. That is the type of analysis pursued in Section~\ref{subsect:satellite.pursuing}, which reveals a severe and pernicious practical manifestation of false confidence.

% , because the other major point established by the false confidence theorem is that there is no theoretical limit on how severely false confidence will manifest itself. Or more precisely, there is no such limit that applies to all measurable propositions. 

% Practical manifestations of false confidence do not involve arbitrarily small sets knowingly drawn around the true value of the parameter being inferred; they involve fixed measurable propositions determined by the structure of the problem itself. 

Additionally, the proof of the false confidence theorem has a direct connection to the practical manifestation of false confidence seen in satellite conjunction analysis. They are both driven by the fact that the size of a measurable set puts an upper bound on the belief that can be assigned to it, which due to additivity, puts a lower bound on the belief that will necessarily be assigned to its complement. The quantitative relationship between epistemic probability and the relative size of the proposition to which it is assigned is perhaps best illustrated by Figure~\ref{fig:ShrinkingProbability}. This relationship holds whether or not the proposition of interest (e.g., collision) includes the true parameter value, but in the proof of the false confidence theorem, it is necessary to show that it does indeed hold, even if $\theta \in B$. That being said, it is important not to over-interpret this commonality.

Practical manifestations of severe false confidence are not limited to problems in which the proposition of interest corresponds to an obviously small set. This is explored in recent and on-going work referenced in Section~\ref{sect:futureWork}. Even in the present work, consider the fact that---while the proposition that the two satellites will collide is Lebesgue-measurable in displacement space---in the original state space, collision would correspond to an infinite hypercylinder, which has an unbounded Lebesgue measure. The false confidence theorem puts a lower bound on how widely severe false confidence will manifest itself in practice, but it does not in any way put an upper bound on the pervasiveness of this phenomenon.

At the end of the day, any risk analyst who plans to represent the results of a statistical inference using an additive belief function must determine whether their proposition of interest (e.g., ``Will these two satellites collide?'') is one of those affected by false confidence. If so, the next question is, ``How severe?'', because there is no general limit to how bad it can get. That is the message of the false confidence theorem.

\subsection{Broader Implications}
\label{subsect:fct.implications}

False confidence is the inevitable result of treating epistemic uncertainty as though it were aleatory variability. Any probability distribution assigns high probability values to large sets. This is appropriate when quantifying aleatory variability, because any realization of a random variable has a high probability of falling in any given set that is large relative to its distribution. Statistical inference is different; a parameter with a fixed value is being inferred from random data. Any proposition about the value of that parameter is either true or false. To paraphrase Nancy Reid and David Cox,\footnote{The exact quote reads, ``[E]ven if an empirical frequency-based view of probability is not used directly as a basis for inference; it is unacceptable if a procedure yielding regions of high probability in the sense of representing uncertain knowledge would, if used repeatedly, give systematically misleading conclusions.'' \cite[p. 295]{reid2015some}} it is a bad inference that treats a false proposition as though it were true, by consistently assigning it high belief values. That is the defect we see in satellite conjunction analysis, and the false confidence theorem establishes that this defect is universal.

% This finding opens a new front in the debate between Bayesian and frequentist schools of thought in statistics. Traditional arguments over epistemic probability have a distinctively philosophical flavor. They focus on issues like the ontological inappropriateness of epistemic probability distributions \cite{Boole_Laws,venn1866chance}, the unjustified use of prior probabilities \cite{Fisher_1930_Inverse_Probability}, and the hypothetical logical consistency of personal belief functions in highly abstract decision-making scenarios \cite{Finetti_Logiques,Savage_Foundations}. In contrast, the false confidence theorem speaks to the practical consequences of treating epistemic uncertainty probabilistically. And on the question of epistemic probability, it appears that the frequentists are right and the Bayesians are wrong. It is not generally safe or sensible or harmless to try to compute the probability of a non-random event. 

This finding opens a new front in the debate between Bayesian and frequentist schools of thought in statistics. Traditional disputes over epistemic probability have focused on seemingly philosophical issues, such as the ontological inappropriateness of epistemic probability distributions \cite{Boole_Laws,Mayo_Book}, the unjustified use of prior probabilities \cite{Fisher_1930_Inverse_Probability}, and the hypothetical logical consistency of personal belief functions in highly abstract decision-making scenarios \cite{Finetti_Logiques,Savage_Foundations}. Despite these disagreements, the statistics community has long enjoyed a truce sustained by results like the Bernstein--von Mises theorem \cite[Ch. 10]{van2000asymptotic}, which indicate that Bayesian and frequentist inferences usually converge with moderate amounts of data.  

% Bernstein-von Mises theorem, which stipulates that Bayesian and frequentist estimators converge with enough data \cite[Ch. 10]{van2000asymptotic}. 

The false confidence theorem undermines that truce, by establishing that the mathematical form in which an inference is expressed can have practical consequences. This finding echoes past criticisms of epistemic probability leveled by advocates of Dempster--Shafer theory, but those past criticisms focus on the structural inability of probability theory to accurately represent incomplete prior knowledge, e.g., \cite[Ch. 3]{Salicone_Text}. The false confidence theorem is much broader in its implications. It applies to all epistemic probability distributions, even those derived from inferences to which the Bernstein--von Mises theorem would also seem to apply.

Simply put, it is not always sensible, nor even harmless, to try to compute the probability of a non-random event. In satellite conjunction analysis, we have a clear real-world example in which the deleterious effects of false confidence are too large and too important to be overlooked. In other applications, there will be propositions similarly affected by false confidence. The question that one must resolve on a case-by-case basis is whether the affected propositions are of practical interest. For now, we focus on identifying an approach to satellite conjunction analysis that is structurally free from false confidence. % It is a problem for which one must check in any Bayesian posterior, fiducial distribution, or confidence distribution. 

% For now, we focus on identifying an approach to satellite conjunction analysis that is structurally free from false confidence. Probability theory is not the only framework for representing epistemic uncertainty. A cornucopia of alternative axiomatic frameworks exist, most of which are non-additive \cite{Klir_GIT_Book}. While the false confidence theorem definitively establishes that probability theory is an inappropriate framework for representing epistemic uncertainty, this does not pose an insurmountable difficulty. We only need to identify a more appropriate framework. 

\section{The Martin--Liu Validity Criterion}
\label{sect:validity}

The false confidence theorem establishes that, to be structurally free from false confidence, it is a necessary condition that a statistical inference be represented non-additively. However, while necessary, non-additivity is not sufficient. It is trivially easy to imagine a non-additive belief function that is just slightly non-additive and therefore suffers from a false confidence problem as bad or almost as bad as that depicted in Section~\ref{sect:falseConfidenceTheorem}. What we need is a criterion, by which to constrain statistical inference, that is sufficient for eliminating false confidence or, more precisely, managing the severity with which false confidence manifests itself. 

The Martin--Liu validity criterion states that the aleatory probability with which a false proposition will be assigned some amount of belief must be no more than one minus that level of belief \cite{martin2016book}. Stated mathematically, 
\begin{equation}
\label{eq:valid}
\textnormal{Pro}_{X;\theta}\!\left( \left\{ x : \textnormal{Bel}_{\Theta;x}\!\left( A \right) \geq 1 - \alpha  \right\}  \right) \leq \alpha , \quad \forall \, \alpha \in \left[ 0 , 1 \right] , \, A \subset \Omega_{\theta} , \, s.t. \,  A \not\ni \theta.
\end{equation}
Since the Martin--Liu validity criterion explicitly limits the rate at which belief is assigned to false propositions, it is almost a tautology that a statistical approach satisfying this criterion will not suffer from the severe false confidence phenomenon described in Sections~\ref{sect:satellite}--\ref{sect:falseConfidenceTheorem}.

This part of the paper implements the Martin--Liu validity criterion as a way to identify statistical methods that are free from false confidence. Section~\ref{subsect:validity.confidence} proves the Martin--Liu validity of confidence intervals and confidence regions. Section~\ref{subsect:validity.ellipsoids} shows that $K\sigma$ uncertainty ellipsoids can often be interpreted as confidence regions. This provides a statistically reliable tool for managing collision risk in satellite conjunction analysis. Section~\ref{subsect:validity.two-dimensional} provides confidence ellipses on the two-dimensional displacement between the two satellites at closest approach.

\subsection{Confidence Regions}
\label{subsect:validity.confidence}

Frequentist confidence intervals and confidence regions have traditionally been rationalized via coverage probability. That is, a $\left( 1 - \alpha \right) \times 100\%$ confidence interval or confidence region is accepted as valid because, over an infinite number of repeated independent draws of the data, at least $\left( 1 - \alpha \right) \times 100\%$ of the resulting confidence intervals or confidence regions would cover the true parameter value. However, if evaluated as tools for assigning belief to propositions, i.e., sets, confidence intervals and confidence regions also satisfy the Martin--Liu validity criterion.

Let $\Gamma_{\alpha}\!\left( x \right)$ be the $\left( 1 - \alpha \right)$ confidence region for $\theta$ given the realized data $x$. The natural assignment of belief and plausibility made by this confidence region are
\[
\textnormal{Bel}_{\Theta;x}\!\left( A \right) = \begin{cases} 1 - \alpha , & A \supset \Gamma_{\alpha}\!\left( x \right) \\ 0 , & A \not\supset \Gamma_{\alpha}\!\left( x \right) \end{cases} \quad \text{and} \quad 
\textnormal{Pls}_{\Theta;x}\!\left( A \right) = \begin{cases} 1 , & A \cap \Gamma_{\alpha}\!\left( x \right) \neq \emptyset \\ \alpha , & A \cap \Gamma_{\alpha}\!\left( x \right) = \emptyset \end{cases}.
\]
That is to say, a confidence region represents the simple assertion that we are $1-\alpha$ confident that the true value of $\theta$ is somewhere inside $\Gamma_{\alpha}\!\left( x \right)$. Any sets containing $\Gamma_{\alpha}\!\left( x \right)$ inherit that confidence; all other sets accrue no positive confidence. It is a coarse way to represent a statistical inference, perhaps more so than is necessary, but it is also safe in the sense that the Martin--Liu validity criterion holds. Indeed, by the coverage probability definition of a confidence region, the event $\{\Gamma_{\alpha}\!\left( x \right) \ni \theta\}$ has $\textnormal{Pro}_{X;\theta}$-probability at least $1-\alpha$. So, for any set $A \not\ni \theta$, the event $\{\Gamma_\alpha(x) \subset A\}$ has $\textnormal{Pro}_{X;\theta}$-probability at most $\alpha$. Therefore, for any false proposition, i.e., any set $A$ such that $A \not\ni \theta$, the probability that said proposition will be assigned a confidence of at least $1-\alpha$ is less than or equal to $\alpha$, hence \eqref{eq:valid} holds. 

When interpreted this way, a confidence region represents a consonant confidence structure, i.e., a possibility distribution with coverage probability properties \cite{Balch_Conf_2012}. A recent paper already proves that all consonant confidence structures satisfy the Martin--Liu validity criterion \cite{denoeux2017frequency}. The reason that we have included the special case proof for confidence regions is that most readers are likely to be unfamiliar with possibility distributions and confidence structures. In contrast, simple confidence intervals and confidence regions are well-established, well-disseminated tools of statistical inference. As proven above, they do not suffer from false confidence.\footnote{It may be helpful to remind the reader that, as pointed out in Section~\ref{sect:falseConfidenceTheorem}, confidence distributions in the sense of \cite{Balch_Conf_2012} do suffer from false confidence, despite also being rationalized in terms of coverage probability. Consonance and coverage probability together are adequate to ensure Martin--Liu validity, but coverage probability alone is not.}

\subsection{Uncertainty Ellipsoids in Satellite Conjunction Analysis}
\label{subsect:validity.ellipsoids}

Uncertainty ellipsoids are already well-established as a screening tool in conjunction analysis \cite{newman2014evolution,helpSTK_threatvolumes}. We propose that they can be used as a principal risk management tool. $K \sigma$ ellipsoids, also known as covariance ellipsoids, can often be treated as approximate confidence regions. The case for this interpretation starts by treating a satellite's true state, $\theta$, as a fixed unknown with $\hat{\theta}\!\left( x \right)$ as an estimate subject to random errors. If a sufficiently large amount of data is used, the error between a satellite's filter-derived state estimate and its true state can often be treated as approximately normal, even if the observation errors are drawn from non-normal distributions \cite{spall1984asymptotic} and even if the filter used is sampling-based \cite{del1999central}. That is, in conjunction analysis, one can often assume
\[
\hat{\theta}\!\left( X \right) - \theta \sim n\!\left( 0 , C_{\Theta} \right)
\]
where $n\!\left( 0 , C_{\Theta} \right)$ is a multivariate normal distribution with zero mean and covariance matrix $C_{\Theta}$. This assumption may break down if a conjunction analysis is done too far in advance \cite{ghrist2012impact} or if too few data are used to estimate the satellite state or if there is an unaccounted for bias error in the filter algorithm or data. For simplicity, we persist in the normality and unbiasedness assumptions. Relaxing these assumptions is left as a topic for future work. 

This representation of aleatory variability in $\hat{\theta}\!\left( X \right) - \theta$ can be eigendecomposed as 
\[
\hat{\theta}\!\left( X \right) - \theta = E_{C_{\Theta}} \Lambda_{C_{\Theta}}^{1/2} \xi
\]
where $X$ is a (hypothetical) random realization drawn from the same distribution as the data used to estimate the orbits, $E_{C_{\Theta}}$ is a matrix whose columns are the eigenvectors of $C_{\Theta}$, $\Lambda_{C_{\Theta}}$ is a diagonal matrix whose entries are the corresponding eigenvalues of $C_{\Theta}$, and $\xi$ is a unit normal vector of the same dimensionality as $\theta$. This can easily be rearranged into a pivot, as
\[
\Lambda_{C_{\Theta}}^{-1/2} E_{C_{\Theta}}^\top \left( \hat{\theta}\!\left( X \right) - \theta \right) = \xi.
\]
%
%where $1/\Lambda_{C_{\Theta}}$ is a diagonal matrix whose entries are the inverse of the eigenvalues of $C_{\Theta}$. 
%Note that $E_{C_{\Theta}}^\top$ and $\sqrt{1/\Lambda_{C_{\Theta}}}$ are simply the inverses of $E_{C_{\Theta}}$ and $\sqrt{\Lambda_{C_{\Theta}}}$, respectively. 
A pivot is a function of the parameter being inferred and the random data used to infer it whose output is a random variable with a known fixed distribution \cite{Casella_Berger_Text}; in this case, $\xi$ is a unit normal vector. Pivots are special because they can be used to derive confidence regions. 

In particular, we are interested $K\sigma$ uncertainty ellipsoids. A $K \sigma$ uncertainty ellipsoid on $\theta$ is generated by projecting a sphere in $\mathbf{\xi}$-space centered at the origin. Let $K > 0$ be some constant. A confidence region of the form
\[
\Gamma_{\alpha}\!\left( x \right) = \left\{ \vartheta : \xi = \Lambda_{C_{\Theta}}^{-1/2} E_{C_{\Theta}}^\top \left( \hat{\theta}\!\left( x \right) - \vartheta \right) , \sum_{i=1}^{\dim(\theta)}{ \xi_i^2 } \leq K^2 \right\}
\]
will map to $\theta$ as a $\dim(\theta)$-dimensional ellipsoid. The confidence, $( 1 - \alpha )$, associated with this uncertainty ellipsoid can be computed as
\[
1 - \alpha = \textnormal{Pro}_{X;\theta}\!\left( \left\{ x : \Gamma_{\alpha}\!\left( x \right) \ni \theta \right\} \right) = \textnormal{Pro}_{\xi}\!\left( \left\{ \xi : \sum^{\dim(\theta)}_{i=1}{ \xi^2_i } \leq K^2 \right\} \right) = F_{\chi^2;\dim(\theta)}\!\left( K^2 \right) 
\]
where $F_{\chi^2;j}$ is the central $\chi_j^2$ distribution function. 
%cumulative distribution function for a central $\chi^2$ distribution with $j$ degrees of freedom. 

Software is commercially available that is capable of computing three-dimensional $K\sigma$ uncertainty ellipsoids representing uncertainty in the predicted position of a satellite and then determining whether or not the resulting ellipsoids for two satellites intersect \cite{helpSTK_threatvolumes,alfano2004eigenvalue}. If two ellipsoids overlap in a conjunction analysis, a collision avoidance maneuver may be needed. If they do not, the two satellites are considered safe, to within the confidence level associated with the uncertainty ellipsoids. 
 
There are a couple of caveats to the approach as currently implemented in the field. First, as traditionally conceptualized, ellipsoid overlap approaches fail to account for the physical size of the two satellites \cite{alfano2003determining}. To correct for this, effective ellipsoid overlap should be defined as occurring when the minimum distance between the two uncertainty ellipsoids is less than the combined radius of the two satellites. So long as the position uncertainties are much larger than the satellites' combined size, the distinction is negligible. However, satellite size should not be ignored in problems with relatively small position uncertainties. 

The other caveat is that, because the two satellites are treated separately, the joint confidence associated with having captured both satellite positions in their respective uncertainty ellipsoids is usually slightly lower than the confidence individually assigned to each ellipsoid. If each ellipsoid has a $1-\alpha$ confidence attached to it, assuming independence between the data used to infer the two satellite trajectories, the confidence attached to simultaneous coverage is
\[
1 - \alpha' = \textnormal{Pro}_{X_1,X_2;\theta_1,\theta_2}\!\left( \left\{ x_1 , x_2 : \Gamma_{\alpha,1}\!\left( x_1 , A \right) \ni \theta_1 , \Gamma_{\alpha,2}\!\left( x_2 , A \right) \ni \theta_2 \right\} \right) = \left( 1 - \alpha \right)^2.
\]
In reality, the errors in the data used to estimate the two satellite positions may not be entirely independent, but the specter of non-independence is easily resolved. Using Fr\'echet bounds \cite{Ferson_Dependency_Sandia}, even under totally unknown dependence between the data,
\[
1 - \alpha' \geq \max\!{ \left( 0 , 1 - 2\alpha \right) }.
\]
Note that if $\alpha$ is small, $\left( 1 - \alpha \right)^2 = 1 - 2\alpha + \alpha^2 \approx 1 - 2\alpha$. So, using $1-2\alpha$ to represent the joint confidence associated with the two uncertainty ellipses is both conservative and extremely close to the independence case. Thanks to the result about confidence regions presented in Section~\ref{subsect:validity.confidence}, so long as one performs a maneuver whenever the two uncertainty ellipsoids intersect, the rate at which collisions occur over a large number of conjunctions---i.e., the operational aleatory probability of collision---will be capped at $\alpha' = 2\alpha$. 

For example, suppose one were to use a pair of $4\sigma$ ellipsoids to represent uncertainty in the positions of the two conjuncting satellites. This represents a $99.9\%$ confidence region around each satellite position. Moreover, one can say with $99.8\%$ confidence that both satellite positions are contained within their respective ellipsoids. And, more importantly, if a satellite operator performs a collision avoidance maneuver whenever those $4\sigma$ ellipsoids overlap, the operational probability of collision will be capped at $0.227\%$.

All of that being said, ellipsoid overlap does not mean that two satellites are definitely going to collide. Rather, it means that collision is still plausible in light of the data. Performing a collision avoidance maneuver is one way of driving down the plausibility of collision, but it is not the only way. For conjunctions predicted several days in advance, prudent navigators will opt to gather more data, which will allow them to shrink the uncertainty ellipsoids before finally deciding whether or not to make a collision avoidance maneuver. So long as they do not wait too long, this is usually the most reasonable course of action. 

\subsection{Uncertainty Ellipses on Displacement}
\label{subsect:validity.two-dimensional}

An alternative approach based on confidence regions is to derive an uncertainty ellipse on the two-dimensional displacement introduced in Section~\ref{subsect:satellite.compute}. Assuming that variability in the direction of the relative velocity vector is small enough to be ignored, the (random) estimate of displacement at closest approach, $\left( \hat{U}' , \hat{V}' \right)$, has a bivariate normal distribution whose mean is equal to the unknown true displacement, $\left( u' , v' \right)$. That is, 
\[
\begin{bmatrix} \hat{U}' \\ \hat{V}' \end{bmatrix} = \begin{bmatrix} u' \\ v' \end{bmatrix} + E_S \begin{bmatrix} S_1 & 0 \\ 0 & S_2 \end{bmatrix} \begin{bmatrix} \xi_1 \\ \xi_2 \end{bmatrix}
\]
where $\xi_1$ and $\xi_2$ are independent unit normal variates and $E_S$, $S_1$, and $S_2$ are the eigendecomposition of the two-dimensional displacement covariance matrix, $C_{\Delta;1:2,1:2}$, defined in Section~\ref{subsect:satellite.compute}. A $K \sigma$ confidence ellipse can be constructed on $\left( u' , v' \right)$ by projecting a circle, $\xi_1^2 + \xi_2^2 \leq K^2$, onto $\left( u' , v' \right)$. The confidence associated with this ellipse is computed from a chi-squared distribution with two degrees of freedom, as outlined in Section~\ref{subsect:validity.ellipsoids}. Recall that collision is associated with the condition $u'^2 + v'^2 \leq R^2$, where $R$ is the combined ``radius'' of the two satellites. So long as the confidence ellipse for $\left( u' , v' \right)$ does not intersect this circle, then one can say with $\left( 1 - \alpha \right)$ confidence that collision will not occur.

The advantage of this approach is that it will yield a lower false alarm rate than the ellispoid overlap approach described in Section~\ref{subsect:validity.ellipsoids}. The disadvantage is that it involves two assumptions that the ellipsoid overlap approach does not, namely that variability in the direction of the estimated relative velocity vector is negligible and that the covariance between the position estimate errors for the two satellites is known. This second assumption is often taken for granted, but it is non-trivial. Additionally, supporting software for the two-dimensional confidence ellipse on displacement is not yet widely available to satellite navigators.
%; although, that may change in time. 

\section{Conclusions}
\label{sect:conclusions}

The work presented in this paper has been done from a fundamentally frequentist point of view, in which $\theta$ (e.g., the satellite states) is treated as having a fixed but unknown value and the data, $x$, (e.g., orbital tracking data) used to infer $\theta$ are modeled as having been generated by a random process (i.e., a process subject to aleatory variability). Someone fully committed to a subjectivist view of uncertainty \cite{Finetti_Logiques,Savage_Foundations} might contest this framing on philosophical grounds. Nevertheless, what we have established, via the false confidence phenomenon, is that the practical distinction between the Bayesian approach to inference and the frequentist approach to inference is not so small as conventional wisdom in the statistics community currently holds. Even when the data are such that results like the Bernstein-von Mises theorem ought to apply, the mathematical form in which an inference is expressed can have large practical consequences that are easily detectable via a frequentist evaluation of the reliability with which belief assignments are made to a proposition of interest (e.g., ``Will these two satellites collide?'').

Our rationale for framing conjunction analysis in frequentist terms is that, as established in the opening of Section~\ref{sect:intro}, everyone in the space industry has an interest in limiting the number of in-orbit collisions. Thus, our goal has been to help satellite operators identify tools adequate for limiting the literal frequency with which collisions involving operational satellites occur. Framing the problem in frequentist terms enables us to do that, whereas framing the problem in Bayesian terms would not. The only circumstance in which a Bayesian analysis could directly enable satellite operators to control the frequency with which collisions occur is if it were based on an aleatory prior \cite{vonMises_1942_BayesFormula} on the satellite states, and this prior would need to reflect the underlying risk of collision due to orbital crowding. As currently practiced, no such aleatory prior is used in conjunction analysis, nor, in our estimation, is one likely to become available in the coming years. An estimate of the aggregate collision risk per unit time seems feasible, but how to parse that into priors on the satellite states is non-obvious. Such an operation may not be well-posed. So, for someone interested in limiting the literal frequency with which collisions occur, it is necessary to treat the satellite states in each conjunction as a fixed, albeit uncertain, reality and then to assess how reliably a proposed risk metric performs. That is the analysis pursued in Section~\ref{subsect:satellite.pursuing}, and under that analysis, epistemic probability of collision does not appear to be a viable risk metric.

There are other engineers and applied scientists tasked with other risk analysis problems for which they, like us, will have practical reasons to take the frequentist view of uncertainty. For those practitioners, the false confidence phenomenon revealed in our work constitutes a serious practical issue. In most practical inference problems, there are uncountably many propositions to which an epistemic probability distribution will consistently accord a high belief value, regardless of whether or not those propositions are true. Any practitioner who intends to represent the results of a statistical inference using an epistemic probability distribution must at least determine whether their proposition of interest is one of those strongly affected by the false confidence phenomenon. If it is, then the practitioner may, like us, wish to pursue an alternative approach. 

Fortunately, the Martin--Liu validity criterion provides a formal benchmark by which to identify statistical methods that are inherently free from severe false confidence. For example, in Section~\ref{subsect:validity.confidence}, we proved that frequentist confidence intervals and confidence regions satisfy the Martin--Liu validity criterion, thus identifying a widely-disseminated and well-established statistical tool that is free from severe false confidence. Technically, this is a special case of a theorem already proved in \cite{denoeux2017frequency}, but our special case proof is offered in the hope that it will be more accessible to a wider readership. In satellite conjunction analysis, this means that $K\sigma$ uncertainty ellipsoids and uncertainty ellipses, which can often be rationalized as confidence regions, are statistically reliable in the same way that epistemic probability of collision is not.

\section{Future and On-Going Work}
\label{sect:futureWork}

Much work remains to be done exploring the range of real-world uncertainty quantification problems that are practically impacted by false confidence. Risk analysis problems with bounded failure domains are strong candidates for severe false confidence. Recent work demonstrates that false confidence can also arise as the result of non-linear uncertainty propagation, even if a marginalization-specific reference posterior is used \cite{carmichael2018exposition,martin2019false}.

As mentioned in Section~\ref{subsect:validity.confidence}, though statistically reliable, confidence intervals and confidence regions are also coarse. More nuanced possibilistic solutions to problems of statistical inference can be derived using consonant confidence structures or Martin--Liu inferential models. In satellite conjunction analysis, the next natural step is to evaluate the statistical performance of the alternative collision risk metrics mentioned in Section~\ref{subsect:intro.modern} to see which, if any, satisfy the Martin--Liu validity criterion.

Finally, a recent conference paper \cite{hejduk2019satellite} asserts that the ellipse overlap approach will motivate too many collision avoidance maneuvers. The authors of \cite{hejduk2019satellite} argue for epistemic probability of collision as a risk metric on the grounds that it is reasonable to treat conjunctions with diluted collision probabilities as akin to the large number of unrecognized conjunctions that occur daily with objects that are too small to be tracked at all, for which no mitigation actions are possible. They view both situations as part of the background risk that is simply accepted as part of operating a satellite. This argument may initially seem attractive, but we see it as promoting an ill-advised tolerance to blind risk. Regardless of the burden incurred by acting on it, the approach we advocate enables satellite operators to quantify and potentially manage their collision risk exposure in the face of high relative uncertainties (i.e.,~high $S/R$ ratios); epistemic probability of collision does not. We recognize the natural desire to balance the goal of preventing collisions against the goal of keeping maneuvers at a reasonable level, and we further recognize that it may not be possible to achieve an acceptable balance between these two goals using present tracking resources. Nevertheless, committing to a risk metric that cannot reliably distinguish ambiguity from confirmed safety is not, in our view, a durable resolution to this impasse.

% Nevertheless, adopting a standard of practice in which conjunctions with high relative uncertainties are effectively ignored is not, in our view, a durable resolution to that impasse.

% , written in part as a response to a review draft of the present work,

% advocated in Section~\ref{sect:validity}

% the standard advocated in \cite{hejduk2019satellite} is to effectively ignore conjunctions with high relative uncertainties, and in our view, that is not a viable long-term resolution to this impasse.

% That being said, we do not see how the standard advocated in \cite{hejduk2019satellite} can possibly be sufficient to prevent or even forestall the day, predicted in \cite{kessler2010kessler}, when accumulated orbital debris will effectively deny humanity access to spaceflight. 

%\ethics{Insert ethics text here.}

\section*{Competing Interest}

Alexandria Validation Consulting, LLC has developed a proprietary algorithm for computing collision risk in satellite conjunction analysis, based in part on insights derived from the research presented in this paper. To protect the financial interests of M.S.B. and his co-investors, this algorithm is currently being withheld from publication.

\section*{Acknowledgments}
Thanks go to the anonymous reviewers for the critical-yet-helpful feedback; to Salvatore Alfano for providing M.S.B. with an introduction to the literature on satellite conjunction analysis and probability dilution; and to  William Oberkampf for his insights and periodic reminders of the systemic threat posed by in-orbit collisions. Special---though, anonymous---thanks go to the satellite navigator who initially brought satellite conjunction analysis to M.S.B.'s attention.
%Thanks go to Salvatore Alfano for providing M.S.B. with an introduction to the literature on satellite conjunction analysis and probability dilution. Thanks also go to William Oberkampf for providing insights on the manuscript and periodic reminders of the systemic threat posed by in-orbit collisions. We would also like to thank Arona Kessler, Julie Balch, and Jonathan Sadeghi for providing their thoughts on earlier versions of the manuscript. Perhaps most importantly, special---though, anonymous---thanks go to the satellite navigator who initially brought satellite conjunction analysis to M.S.B.'s attention.

\bibliography{generalRef}

\begin{thebibliography}{10}
\expandafter\ifx\csname url\endcsname\relax
  \def\url#1{\texttt{#1}}\fi
\expandafter\ifx\csname urlprefix\endcsname\relax\def\urlprefix{URL }\fi
\expandafter\ifx\csname href\endcsname\relax
  \def\href#1#2{#2} \def\path#1{#1}\fi

\bibitem{kessler2010kessler}
D.~J. Kessler, N.~L. Johnson, J.~Liou, M.~Matney, The {Kessler Syndrome}:
  implications to future space operations, {AAS}-10-016, in: 33rd Annual {AAS}
  Guidance and Control Conference, American Astronautical Society,
  Breckenridge, CO, 2010.

\bibitem{liou2006risks}
J.-C. Liou, N.~L. Johnson, Risks in space from orbiting debris, Science 311
  (2006) 340--341.

\bibitem{liou2008instability}
J.-C. Liou, N.~L. Johnson, Instability of the present leo satellite
  populations, Advances in Space Research 41~(7) (2008) 1046--1053.

\bibitem{alfano2003collision}
S.~Alfano, Relating position uncertainty to maximum conjunction probability,
  {AAS}-03-548, in: AAS/AIAA Astrodynamics Specialists Conference, Big Sky, MT,
  2003.

\bibitem{frigm2009single}
C.~Frigm, A single conjunction risk assessment metric: the f-value, Advances in
  the Astronautical Sciences 135~(2) (2009) 1175--1192.

\bibitem{plakalovic2011tuned}
D.~Plakalovi{\'c}, M.~Hejduk, R.~Frigm, L.~Newman, A tuned single parameter for
  representing conjunction risk, {AAS}-11-430, in: 2011 AIAA/AAS Astrodynamics
  Specialist Conference, Girdwood, AK, 2011.

\bibitem{sun2014spacecraft}
Z.~Sun, Y.~Luo, Z.~Niu, Spacecraft rendezvous trajectory safety quantitative
  performance index eliminating probability dilution, Science China
  Technological Sciences 57~(6) (2014) 1219--1228.

\bibitem{balch2016corrector}
M.~S. Balch, A corrector for probability dilution in satellite conjunction
  analysis, {AIAA}-2016-1445, in: 18th AIAA Non-Deterministic Approaches
  Conference, San Diego, CA, 2016.

\bibitem{carpenter2017relevance}
J.~R. Carpenter, S.~Alfano, D.~T. Hall, M.~D. Hejduk, J.~A. Gaebler, M.~K. Jah,
  S.~O. Hasan, R.~L. Besser, R.~R. DeHart, M.~G. Duncan, M.~S. Herron, W.~J.
  Guit, Relevance of the american statistical society's warning on p-values for
  conjunction assessment, {AAS}-17-614, in: AAS/AIAA Astrodynamics Specialist
  Conference, Stevenson, WA, 2017.

\bibitem{newman2014evolution}
L.~K. Newman, R.~C. Frigm, M.~G. Duncan, M.~D. Hejduk, Evolution and
  implementation of the nasa robotic conjunction assessment risk analysis
  concept of operations, in: Advanced Maui Optical and Space Surveillance
  Technologies Conference, Wailea, HI, 2014.

\bibitem{Barnett_Comparative}
V.~Barnett, Comparative Statistical Inference, 3rd Edition, Wiley and Sons,
  Chichester, NY, 1999.

\bibitem{Laplace_Essay}
P.~S. Laplace, A Philosophical Essay on Probabilities, Wiley and Sons, New
  York, NY, Translated 1902; Original 1814.

\bibitem{Finetti_Logiques}
B.~de~Finetti, La pr\'evision: ses lois logiques, ses sources subjectives,
  Annales de l'Institute Henri Poincar$\acute{e}$, Paris, France, 1935.

\bibitem{robert2007bayesian}
C.~P. Robert, The Bayesian Choice: from Decision-Theoretic Foundations to
  Computational Implementation, Springer Science \& Business Media, New York,
  NY, 2007.

\bibitem{Boole_Laws}
G.~Boole, An Investigation of the Laws of Thought, Project Gutenberg,
  www.gutenberg.net, Online release 2005; Original 1854.

\bibitem{neyman1937outline}
J.~Neyman, Outline of a theory of statistical estimation based on the classical
  theory of probability, Philosophical Transactions of the Royal Society of
  London. Series A, Mathematical and Physical Sciences 236~(767) (1937)
  333--380.

\bibitem{Mayo_Book}
D.~G. Mayo, Error and the growth of experimental knowledge, University of
  Chicago Press, Chicago, IL, 1996.

\bibitem{ferson1996different}
S.~Ferson, L.~R. Ginzburg, Different methods are needed to propagate ignorance
  and variability, Reliability Engineering \& System Safety 54~(2-3) (1996)
  133--144.

\bibitem{Salicone_Text}
S.~Salicone, Measurement Uncertainty: An Approach via the Mathematical Theory
  of Evidence, Springer Science \& Business Media, LLC, New York, New York,
  2007.

\bibitem{Roy_Oberkampf_Book}
W.~L. Oberkampf, C.~J. Roy, Verification and Validation in Scientific
  Computing, Cambridge University Press, New York, NY, 2010.

\bibitem{Shafer_1976_Book}
G.~Shafer, A Mathematical Theory of Evidence, Prinecton University Press,
  Princeton, NJ, 1976.

\bibitem{Balch_Conf_2012}
M.~S. Balch, Mathematical foundations for a theory of confidence structures,
  International Journal of Approximate Reasoning 53 (2012) 1003--1019.

\bibitem{Schweder_Hjort_book}
T.~Schweder, N.~L. Hjort, Confidence, Likelihood, Probability: Statistical
  Inference with Confidence Distributions, Cambridge University Press, New
  York, NY, 2016.

\bibitem{Singh_2013_CD}
M.~Xie, K.~Singh, Confidence distribution, the frequentist distribution
  estimator of a parameter: a review, Int. Stat. Rev. 81~(1) (2013) 3--39.

\bibitem{vonMises_1942_BayesFormula}
R.~von Mises, On the correct use of bayes' formula, The Annals of Mathematical
  Statistics 13~(2) (1942) 156 -- 165.

\bibitem{martin2016book}
R.~Martin, C.~Liu, Inferential Models: Reasoning with Uncertainty, CRC Press,
  Boca Raton, FL, 2016.

\bibitem{alfano2016probability}
S.~Alfano, D.~Oltrogge, Probability of collision: Valuation, variability,
  visualization and validity, {AIAA}-2016-5654, in: AIAA/AAS Astrodynamics
  Specialist Conference, Long Beach, CA, 2016.

\bibitem{frigm2015total}
R.~C. Frigm, M.~D. Hejduk, L.~C. Johnson, D.~Plakalovic, Total probability of
  collision as a metric for finite conjunction assessment and collision risk
  management, in: Proceedings of the Advanced Maui Optical and Space
  Surveillance Technologies Conference, Wailea, HI, 2015.

\bibitem{alfriend1999probability}
K.~T. Alfriend, M.~R. Akella, J.~Frisbee, J.~L. Foster, D.-J. Lee, M.~Wilkins,
  Probability of collision error analysis, Space Debris 1~(1) (1999) 21--35.

\bibitem{patera2001general}
R.~P. Patera, General method for calculating satellite collision probability,
  Journal of Guidance, Control, and Dynamics 24~(4) (2001) 716--722.

\bibitem{chan2008spacecraft}
F.~K. Chan, Spacecraft collision probability, Aerospace Press, El Segundo, CA,
  2008.

\bibitem{meinhold1983understanding}
R.~J. Meinhold, N.~D. Singpurwalla, Understanding the kalman filter, The
  American Statistician 37~(2) (1983) 123--127.

\bibitem{hannig2009generalized}
J.~Hannig, On generalized fiducial inference, Statistica Sinica 19 (2009)
  491--544.

\bibitem{hall2017time}
D.~T. Hall, M.~D. Hejduk, L.~C. Johnson, Time dependence of collision
  probabilities during satellite conjunctions, {AAS}-17-271, in: {AAS/AIAA}
  Space Flight Mechanics Meeting, San Antonio, TX, 2017.

\bibitem{helpSTK_nonLinearProb}
{Analytical Graphics, Inc.}, {STK}: Nonlinar probability tool, accessed online
  at http://help.agi.com/stk/index.htm\#cat/Cat03-08.htm (2018).

\bibitem{patera2005calculating}
R.~P. Patera, Calculating collision probability for arbitrary space-vehicle
  shapes via numerical quadrature, Journal of Guidance, Control, and Dynamics
  28~(6) (2005) 1326--1331.

\bibitem{sabol2010linearized}
C.~Sabol, T.~Sukut, K.~Hill, K.~T. Alfriend, B.~Wright, Y.~Li, P.~Schumacher,
  Linearized orbit covariance generation and propagation analysis via simple
  monte carlo simulations, Tech. rep., Texas Engineering Experiment Station,
  College Station, TX (2010).

\bibitem{ghrist2012impact}
R.~W. Ghrist, D.~Plakalovic, Impact of non-gaussian error volumes on
  conjunction assessment risk analysis, {AIAA}-2012-4965, in: AIAA/AAS
  Astrodynamics Specialist Conference, Minneapolis, MN, 2012.

\bibitem{newman2010overview}
L.~K. Newman, The nasa robotic conjunction assessment process: Overview and
  operational experiences, Acta Astronautica 66~(7) (2010) 1253--1261.

\bibitem{hannig2016generalized}
J.~Hannig, H.~Iyer, R.~C. Lai, T.~C. Lee, Generalized fiducial inference: A
  review and new results, Journal of the American Statistical Association
  111~(515) (2016) 1346--1361.

\bibitem{Fisher_1935_Fiducial_Argument}
R.~A. Fisher, The fiducial argument in statistical inference, Annals of
  Eugenics 6 (1935) 391--398.

\bibitem{reid2015some}
N.~Reid, D.~R. Cox, On some principles of statistical inference, International
  Statistical Review 83~(2) (2015) 293--308.

\bibitem{Fisher_1930_Inverse_Probability}
R.~A. Fisher, Inverse probability, Proceedings of the Cambridge Philosophical
  Society 26 (1930) 528--535.

\bibitem{Savage_Foundations}
L.~J. Savage, The Foundations of Statistics, 2nd Edition, Dover Publications
  Inc., New York, NY, 1972.

\bibitem{van2000asymptotic}
A.~W. Van~der Vaart, Asymptotic statistics, Cambridge University Press, New
  York, NY, 2000.

\bibitem{denoeux2017frequency}
T.~Denoeux, S.~Li, Frequency-calibrated belief functions: Review and new
  insights, International Journal of Approximate Reasoning 92 (2018) 232--254.

\bibitem{helpSTK_threatvolumes}
{Analytical Graphics, Inc.}, {STK}: Defining threat volume, accessed online at
  http://help.agi.com/stk/index.htm\#cat/Cat03-01.htm (2018).

\bibitem{spall1984asymptotic}
J.~C. Spall, K.~D. Wall, Asymptotic distribution theory for the kalman filter
  state estimator, Communications in Statistics-Theory and Methods 13~(16)
  (1984) 1981--2003.

\bibitem{del1999central}
P.~Del~Moral, A.~Guionnet, Central limit theorem for nonlinear filtering and
  interacting particle systems, Annals of Applied Probability 9~(2) (1999)
  275--297.

\bibitem{Casella_Berger_Text}
G.~Casella, R.~L. Berger, Statistical Inference, 2nd Edition, Thomson Learning,
  Pacific Grove, CA, 2002.

\bibitem{alfano2004eigenvalue}
S.~Alfano, F.~K. Chan, M.~L. Greer, Eigenvalue quadric surface method for
  determining when two ellipsoids share common volume for use in spatial
  collision detection and avoidance, uS Patent 6,694,283 (Feb.~17 2004).

\bibitem{alfano2003determining}
S.~Alfano, M.~L. Greer, Determining if two solid ellipsoids intersect, Journal
  of Guidance, Control, and Dynamics 26~(1) (2003) 106--110.

\bibitem{Ferson_Dependency_Sandia}
S.~Ferson, R.~B. Nelsen, J.~Hajagos, D.~J. Berleant, J.~Zhang, W.~T. Tucker,
  L.~R. Ginzburg, W.~L. Oberkampf, Dependence in probabilistic modeling,
  dempster-shafer theory, and probability bounds analysis, {SAND}-2004-3072,
  Tech. rep., Sandia National Laboratories, Albuquerque, NM (2004).

\bibitem{carmichael2018exposition}
I.~Carmichael, J.~Williams, An exposition of the false confidence theorem, Stat
  7~(1) (2018) e201.

\bibitem{martin2019false}
R.~Martin, False confidence, non-additive beliefs, and valid statistical
  inference, arXiv preprint arXiv:1607.05051.

\bibitem{hejduk2019satellite}
M.~D. Hejduk, D.~E. Snow, L.~K. Newman, Satellite conjunction assessment risk
  analysis for ``dilution region'' events: Issues and operational approaches,
  in: Space Traffic Management Conference, Austin, TX, 2019.

\end{thebibliography}

\end{document}